\documentclass{amsart} 
\usepackage{amssymb,latexsym,ifthen}

\newcommand{\circled}[1]{{\,\bigcirc\hspace{-.123in} #1\,}}

\def\ii{\mathbf{i}}
\def\jj{\mathbf{j}}
\def\circi{{\bigcirc\!\!\!\! i\,}}
\def\circj{{\bigcirc\!\!\!\! j\,}}
\def\A{\mathcal{A}}
\def\l{\ell}

\newcommand{\light}[1]{{#1}} 
\newcommand{\dark}[1]{{#1}} 
\newcommand{\bulletcolor}[1]{{#1}} 

\newcounter{colorversion}

\ifthenelse{\thecolorversion=1}
{
\usepackage{color}
\definecolor{light}{rgb}{1,0.1,0.1}        
\definecolor{dark}{cmyk}{1,0.4,0,0}  
\definecolor{bulletcolor}{cmyk}{1,0.3,1,0}  
\renewcommand{\light}[1]{\color{light}{#1}\color{black}}
\renewcommand{\dark}[1]{\color{dark}{#1}\color{black}}
\renewcommand{\bulletcolor}[1]{\color{bulletcolor}{#1}\color{black}}
}{}

\newtheorem{theorem}{Theorem}
\newtheorem{proposition}[theorem]{Proposition}
\newtheorem{corollary}[theorem]{Corollary}

\newtheorem{lemma}[theorem]{Lemma}
\newtheorem{example}[theorem]{Example}
\newtheorem{definition}[theorem]{Definition}
\newtheorem{conjecture}[theorem]{Conjecture}

\def\proof{\smallskip\noindent {\bf Proof. }}
\def\endproof{\hfill$\square$\medskip}
\def\endproofmath{\quad\square}

\def\CC{\mathbb{C}}

\def\wnot{w_\mathrm{o}}
\def\ii{\mathbf{i}}

\newcommand{\matbr}[4]{\left[\!\!\begin{array}{cc}
#1 & #2 \\
#3 & #4 \\
\end{array}\!\!\right]}

\begin{document}

\title[Total positivity: tests and parametrizations]
{{\ }\\[-.9in]
Total positivity:\\ tests and parametrizations}

\author{Sergey Fomin}
\address{Department of Mathematics, Massachusetts Institute of
  Technology, Cambridge, Massachusetts 02139}
\email{fomin@math.mit.edu}

\author{Andrei Zelevinsky}
\address{Department of Mathematics, Northeastern University, 
Boston, MA 02115, USA}
\email{andrei@neu.edu}

\date{
August 1999}

\maketitle

\section*{Introduction}

A matrix is \emph{totally positive} 
\ (resp.\ \emph{totally nonnegative})
if all its minors are positive (resp.\ nonnegative) 
real numbers. 
The first systematic study of these classes of matrices was undertaken
in the 1930s by F.~R.~Gantmacher and M.~G.~Krein \cite{GK35, GK37,
  GK}, 
who established their remarkable spectral properties 
(in particular, an $n \times n$ 
totally positive matrix~$x$  
has  $n$ distinct positive eigenvalues).
Earlier, I.~J.~Schoenberg~\cite{schoenberg} discovered the connection
between total nonnegativity and the following \emph{variation-diminishing
  property:}  the number of sign changes in
a vector does not increase upon multiplying by~$x$. 

Total positivity found numerous applications
and was studied from many different angles.
An incomplete list includes:  
oscillations in mechanical systems
(the original motivation in~\cite{GK}), 
stochastic processes and approximation theory 
\cite{gasca-micchelli, karlin}, 
P\'olya frequency sequences \cite{karlin, schoenberg-werke}, 
representation theory of the infinite symmetric group and the
Edrei-Thoma theorem~\cite{edrei, thoma}, 
planar resistor networks~\cite{CIM},
unimodality and log-concavity~\cite{stanley-log}, 
and theory of immanants~\cite{stembridge-immanants}. 
Further references can be found in S.~Karlin's book~\cite{karlin}
and in the surveys~\cite{ando, brenti-survey, pinkus}.

In this paper, we focus on the following two problems: 
\begin{itemize}
\item[(i)] \emph{parametrizing} all totally nonnegative matrices; 

\item[(ii)] \emph{testing} a matrix for total positivity. 
\end{itemize}

\noindent
Our interest in these problems stemmed from a surprising
represent\-ation-theoretic connection between total positivity and
canonical bases for quantum groups, discovered by
G.~Lusztig~\cite{lusztig} (cf.\ also the surveys 
\cite{littelmann, lusztig-survey}). 
Among other things, he extended the subject by defining
totally positive and totally nonnegative elements 
for any reductive group. 
Further development of these ideas in~\cite{BFZ, BZ, FZ, FZosc} aims
at generalizing the whole body of classical determinantal calculus 
to any semisimple group. 

As it often happens, putting things in a more general perspective
shed new light on this classical subject.
In the next two sections 
of this paper, we provide self-contained proofs (many of them new) of
the fundamental results on problems (i)--(ii), 
due to A.~Whitney~\cite{whitney}, C.~Loewner~\cite{loewner}, 
C.~Cryer~\cite{cryer, cryer76}, and M.~Gasca and
J.~M.~Pe\~na~\cite{gasca-pena}. 
The rest of the paper presents more recent results obtained
in~\cite{FZ}: 
a family of efficient total positivity criteria,
and explicit formulas for expanding a generic matrix
into a product of elementary Jacobi matrices.
These results and their proofs can be generalized to arbitrary semisimple
groups~\cite{BZ,FZ}, but we do not discuss this here.

Our approach to the subject relies on two combinatorial 
constructions. 
The first one is well known: it associates a totally nonnegative
matrix to a planar directed graph with positively weighted
edges 
(in fact, every totally nonnegative matrix can be obtained in this
way~\cite{brenti}).   
Our second combinatorial tool was introduced in \cite{FZ}; 
it is a particular class of colored pseudoline arrangements that we
call the \emph{double wiring diagrams.}





\section*{Planar networks}
\label{sec:planar networks}

To the uninitiated, it might be unclear that totally positive
matrices of arbitrary order exist at all. 
As a warm-up, we invite the reader to check that 
every matrix given by
\begin{equation}
\label{eq:3by3}
\left[\begin{array}{ccc}
d    & dh         & dhi \\
bd   & bdh\!+\!e      & bdhi\!+\!eg\!+\!ei \\
abd  & abdh\!+\!ae\!+\!ce & abdhi\!+\!(a\!+\!c)e(g\!+\!i) \!+\! f\\
\end{array}\right] ,
\end{equation}
where the numbers $a,b,c,d,e,f,g,h,i$ are positive, is totally
positive.
It will follow from the results below that  \emph{every} 
$3 \times 3$ totally positive matrix has this form.

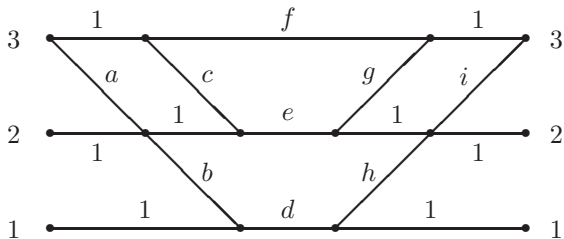
\begin{figure}[ht]
\setlength{\unitlength}{1.8pt} 

\begin{center}
\begin{picture}(100,45)(0,0)
\thicklines

\put(0,0){\line(1,0){100}}
\put(0,20){\line(1,0){100}}
\put(0,40){\line(1,0){100}}

\put(60,20){\line(1,1){20}}
\put(60,0){\line(1,1){40}}

\put(0,40){\line(1,-1){40}}
\put(20,40){\line(1,-1){20}}

  \put(105,-2){${1}$}
  \put(105,18){${2}$}
  \put(105,38){${3}$}

  \put(-9,-2){${1}$}
  \put(-9,18){${2}$}
  \put(-9,38){${3}$}

  \put(0,0){\circle*{1.5}}
  \put(0,20){\circle*{1.5}}
  \put(0,40){\circle*{1.5}}
  \put(20,20){\circle*{1.5}}
  \put(20,40){\circle*{1.5}}
  \put(40,0){\circle*{1.5}}
  \put(40,20){\circle*{1.5}}
  \put(60,0){\circle*{1.5}}
  \put(60,20){\circle*{1.5}}
  \put(80,20){\circle*{1.5}}
  \put(80,40){\circle*{1.5}}
  \put(100,0){\circle*{1.5}}
  \put(100,20){\circle*{1.5}}
  \put(100,40){\circle*{1.5}}

  \put(10,44){\makebox(0,0){$1$}}
  \put(50,44){\makebox(0,0){$f$}}
  \put(90,44){\makebox(0,0){$1$}}
  \put(20, 4){\makebox(0,0){$1$}}
  \put(50, 4){\makebox(0,0){$d$}}
  \put(80, 4){\makebox(0,0){$1$}}
  \put(10,16){\makebox(0,0){$1$}}
  \put(27,24){\makebox(0,0){$1$}}
  \put(50,24){\makebox(0,0){$e$}}
  \put(73,24){\makebox(0,0){$1$}}
  \put(90,16){\makebox(0,0){$1$}}
  \put(33,12){\makebox(0,0){$b$}}
  \put(33,32){\makebox(0,0){$c$}}
  \put(13,32){\makebox(0,0){$a$}}
  \put(67,12){\makebox(0,0){$h$}}
  \put(67,32){\makebox(0,0){$g$}}
  \put(87,32){\makebox(0,0){$i$}}

\end{picture}
\end{center}

\caption{A planar network}
\label{fig:planar}
\end{figure}

\begin{figure}[ht]
\setlength{\unitlength}{1.2pt} 
\begin{center}
\begin{picture}(180,85)(0,-5)
\thicklines

\put(0,0){\line(1,0){80}}
\put(0,20){\line(1,0){80}}
\put(0,40){\line(1,0){80}}
\put(0,60){\line(1,0){80}}
\put(0,80){\line(1,0){80}}

\bulletcolor{
\put(80,0){\line(1,0){20}}
\put(80,20){\line(1,0){20}}
\put(80,40){\line(1,0){20}}
\put(80,60){\line(1,0){20}}
\put(80,80){\line(1,0){20}}
}

\put(100,0){\line(1,0){80}}
\put(100,20){\line(1,0){80}}
\put(100,40){\line(1,0){80}}
\put(100,60){\line(1,0){80}}
\put(100,80){\line(1,0){80}}

\dark{
\put(100,0){\line(1,1){80}}
\put(100,20){\line(1,1){60}}
\put(100,40){\line(1,1){40}}
\put(100,60){\line(1,1){20}} 
}

\thinlines
\light{
\put(60,80){\line(1,-1){20}}
\put(40,80){\line(1,-1){40}}
\put(20,80){\line(1,-1){60}}
\put(0,80){\line(1,-1){80}}
}

\thicklines

\qbezier[40](0,3)(50,3)(99,3)
\qbezier[40](0,23)(50,23)(99,23)
\qbezier[25](101,5)(120,24)(138,42)
\qbezier[25](101,25)(120,44)(138,62)
\qbezier[20](140,43)(160,43)(180,43)
\qbezier[20](140,63)(160,63)(180,63)

\thicklines

\put(115,50){\vector(1,0){10}}

  \put(185,-2){${1}$}
  \put(185,18){${2}$}
  \put(185,38){${3}$}
  \put(185,58){${\vdots}$}
  \put(185,78){$n$}

  \put(-9,-2){${1}$}
  \put(-9,18){${2}$}
  \put(-9,38){${3}$}
  \put(-9,58){${\vdots}$}
  \put(-9,78){$n$}

  \put(0,0){\circle*{2.5}}
  \put(0,20){\circle*{2.5}}
  \put(0,40){\circle*{2.5}}
  \put(0,60){\circle*{2.5}}
  \put(0,80){\circle*{2.5}}

  \put(20,60){\circle*{2.5}}
  \put(20,80){\circle*{2.5}}
  \put(40,40){\circle*{2.5}}
  \put(40,60){\circle*{2.5}}
  \put(40,80){\circle*{2.5}}
  \put(60,20){\circle*{2.5}}
  \put(60,40){\circle*{2.5}}
  \put(60,60){\circle*{2.5}}
  \put(60,80){\circle*{2.5}}
  \put(80,0){\circle*{2.5}}
  \put(80,20){\circle*{2.5}}
  \put(80,40){\circle*{2.5}}
  \put(80,60){\circle*{2.5}}
  \put(80,80){\circle*{2.5}}
  \put(100,0){\circle*{2.5}}
  \put(100,20){\circle*{2.5}}
  \put(100,40){\circle*{2.5}}
  \put(100,60){\circle*{2.5}}
  \put(100,80){\circle*{2.5}}
  \put(120,20){\circle*{2.5}}
  \put(120,40){\circle*{2.5}}
  \put(120,60){\circle*{2.5}}
  \put(120,80){\circle*{2.5}}
  \put(140,40){\circle*{2.5}}
  \put(140,60){\circle*{2.5}}
  \put(140,80){\circle*{2.5}}
  \put(160,60){\circle*{2.5}}
  \put(160,80){\circle*{2.5}}
  \put(180,0){\circle*{2.5}}
  \put(180,20){\circle*{2.5}}
  \put(180,40){\circle*{2.5}}
  \put(180,60){\circle*{2.5}}
  \put(180,80){\circle*{2.5}}

\end{picture}
\end{center}

\caption{Planar network~$\Gamma_0$}
\label{fig:maxnetwork}
\end{figure}
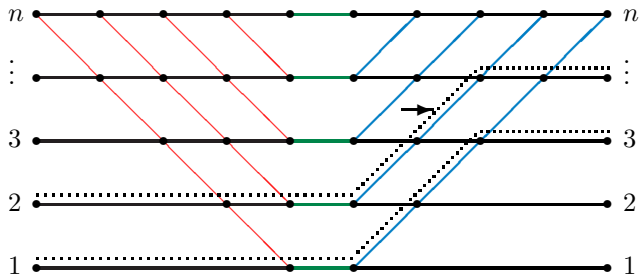

\begin{figure}[ht]
\setlength{\unitlength}{0.8pt} 

\begin{center}
\begin{picture}(80,80)(0,3)
\thicklines

  \put(0,0){\line(1,0){80}}
  \put(0,0){\line(0,1){80}}
  \put(80,0){\line(0,1){80}}
  \put(0,80){\line(1,0){80}}

\dark{
  \put(40,50){\line(1,0){28}}
  \put(40,50){\line(0,1){28}}
  \put(40,78){\line(1,0){28}}
  \put(68,50){\line(0,1){28}}
}

\light{
  \put( 2,10){\line(1,0){58}}
  \put( 2,10){\line(0,1){58}}
  \put( 2,68){\line(1,0){58}}
  \put(60,10){\line(0,1){58}}
}
\end{picture}
\end{center}
\caption{
Initial minors
}
\label{fig:gasca}

\end{figure}
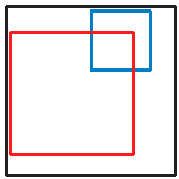

We will now describe a general procedure that
produces totally nonnegative matrices. 
In what follows, a \emph{planar network}~$(\Gamma,\omega)$ is an acyclic
directed planar graph~$\Gamma$ whose edges~$e$ are assigned scalar 
\emph{weights}~$\omega(e)$. 
In all of our examples (Figures~\ref{fig:planar}--\ref{fig:planar-network}), 
we assume the edges of~$\Gamma$ \emph{directed left to right.} 
Also, each of our networks will have $n$ \emph{sources} and $n$
\emph{sinks,} located at the left (resp.\ right) edge of the picture,
and numbered bottom-to-top.

The weight of a directed path in $\Gamma$ is defined as the product
of the weights of its edges. 
The \emph{weight matrix} $x(\Gamma,\omega)$ 
is an $n \times n$ matrix whose $(i,j)$-entry is the sum of weights of
all paths from the source~$i$ to the sink~$j$.
For example, the weight matrix of the network in
Figure~\ref{fig:planar} is given by~(\ref{eq:3by3}). 
  
The minors of the weight matrix of a planar network have an important 
combinatorial interpretation, which can be traced to
B.~Lindstr\"om~\cite{lindstrom} and further to S.~Karlin and
G.~McGregor~\cite{karlin-mcgregor} (implicit), and whose many 
applications were given by I.~Gessel and G.~X.~Viennot~\cite{GV, GV2}. 

In what follows, $\Delta_{I,J}(x)$ denotes the minor of a matrix $x$
with the row set $I$ and the column set~$J$.

The weight of a collection of directed paths in~$\Gamma$ is
defined to be the product of their weights. 

\begin{lemma} {\rm (Lindstr\"om's Lemma)}
\label{lem:planar-minors}
A minor $\Delta_{I,J}$ of the weight matrix 
of a planar network is equal to the sum of
weights of all collections of vertex-disjoint paths 
that connect the sources labeled by~$I$ with the sinks labeled by~$J$.  
\end{lemma}

To illustrate, consider the matrix~$x$ in~(\ref{eq:3by3}). 
We have, e.g., $\Delta_{23,23}(x)=bcdegh+bdfh+fe$, which also equals
the sum of the weights of the three vertex-disjoint path collections in
Figure~\ref{fig:planar} that 
connect sources~2 and~3 to sinks~2 and~3. 

\proof
It suffices to prove the lemma for the determinant of the whole weight
matrix $x=x(\Gamma,\omega)$, i.e., for the case $I = J = [1,n]$. 
Expanding the determinant, we obtain
\begin{equation}
\label{eq:det network}
\det (x) = \sum_w \sum_\pi {\rm sgn} (w)\, \omega (\pi) \ ,
\end{equation}
the sum being over all permutations $w$ in the symmetric
group~$\mathcal{S}_n\,$, 
and over all collections of paths $\pi = (\pi_1, \dots, \pi_n)$
such that $\pi_i$ joins the source $i$ with the sink~$w(i)$.
Any collection $\pi$ of vertex-disjoint paths is
associated with the identity permutation; 
hence $\omega(\pi)$ appears in (\ref{eq:det network}) with the
positive sign.  
We need to show that all other terms in (\ref{eq:det network}) cancel out. 
Deforming $\Gamma$ a bit if necessary, we may assume that no
two vertices lie on the same vertical line. 
This makes the following involution on the non-ver\-tex-disjoint
collections of paths well-defined: take the rightmost point of
intersection of two paths 
in $\pi$, and switch the parts of these paths lying to the right of 
this point. 
This involution preserves the weight of~$\pi$, while changing the sign
of the associated permutation~$w$;
the corresponding pairing of terms in (\ref{eq:det network}) provides
the desired cancellation. 
\endproof

\begin{corollary}
\label{cor:positive planar-minors1}
If a planar network has nonnegative real weights, 
then its weight matrix is totally nonnegative.
\end{corollary}

An aside: note that the weight matrix of the network 
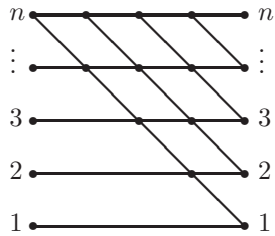
\begin{center}
\begin{picture}(80,92)(0,-5)
\thicklines

\put(0,0){\line(1,0){80}}
\put(0,20){\line(1,0){80}}
\put(0,40){\line(1,0){80}}
\put(0,60){\line(1,0){80}}
\put(0,80){\line(1,0){80}}

\put(60,80){\line(1,-1){20}}
\put(40,80){\line(1,-1){40}}
\put(20,80){\line(1,-1){60}}
\put(0,80){\line(1,-1){80}}

  \put(85,-2){${1}$}
  \put(85,18){${2}$}
  \put(85,38){${3}$}
  \put(85,58){${\vdots}$}
  \put(85,78){$n$}

  \put(-9,-2){${1}$}
  \put(-9,18){${2}$}
  \put(-9,38){${3}$}
  \put(-9,58){${\vdots}$}
  \put(-9,78){$n$}

  \put(0,0){\circle*{2.5}}
  \put(0,20){\circle*{2.5}}
  \put(0,40){\circle*{2.5}}
  \put(0,60){\circle*{2.5}}
  \put(0,80){\circle*{2.5}}

  \put(20,60){\circle*{2.5}}
  \put(20,80){\circle*{2.5}}
  \put(40,40){\circle*{2.5}}
  \put(40,60){\circle*{2.5}}
  \put(40,80){\circle*{2.5}}
  \put(60,20){\circle*{2.5}}
  \put(60,40){\circle*{2.5}}
  \put(60,60){\circle*{2.5}}
  \put(60,80){\circle*{2.5}}
  \put(80,0){\circle*{2.5}}
  \put(80,20){\circle*{2.5}}
  \put(80,40){\circle*{2.5}}
  \put(80,60){\circle*{2.5}}
  \put(80,80){\circle*{2.5}}

\end{picture}
\end{center}
(with unit edge weights) is the ``Pascal triangle'' 
\[
\left[
\begin{array}{cccccc}
1 & 0 & 0 & 0 & 0 & \cdots \\
1 & 1 & 0 & 0 & 0 & \cdots \\
1 & 2 & 1 & 0 & 0 & \cdots \\
1 & 3 & 3 & 1 & 0 & \cdots \\
1 & 4 & 6 & 4 & 1 & \cdots \\
\cdots & \cdots & \cdots & \cdots & \cdots & \cdots \\
\end{array}
\right]\,,
\]
which is totally nonnegative by Corollary~\ref{cor:positive
planar-minors1}. Similar arguments can be used to show total
nonnegativity of various other combinatorial matrices, such as the matrices
of $q$-binomial coefficients, Stirling numbers of both kinds, etc. 

We call a planar network $\Gamma$ \emph{totally connected} if, for any two
subsets $I, J \subset [1,n]$ of  
the same cardinality, there exists a collection of vertex-disjoint 
paths in $\Gamma$ connecting the sources labeled by~$I$ with 
the sinks labeled by~$J$.

\begin{corollary}
\label{cor:positive planar-minors2}
If a totally connected planar network has positive weights,
then its weight matrix is totally positive. 
\end{corollary}

For any $n$, let $\Gamma_0$ denote the network 
shown in Figure~\ref{fig:maxnetwork}. 
Direct inspection shows that $\Gamma_0$ is totally connected. 

\begin{corollary}
\label{cor:minimal network}
For any choice of positive weights $\omega(e)$,
the weight matrix $x(\Gamma_0,\omega)$ 
is totally positive. 
\end{corollary}

It turns out that this construction produces 
\emph{all} totally positive
matrices; this result is essentially equivalent to A.~Whitney's 
Reduction Theorem~\cite{whitney}, and can be sharpened as follows. 
Call an edge of $\Gamma_0$ \emph{essential} if it either is 
slanted, or is one of the $n$ horizontal edges in the middle of the
network. 
Note that  $\Gamma_0$ has exactly $n^2$ 
essential edges. 
A weighting $\omega$ of $\Gamma_0$ is \emph{essential} if $\omega (e) \neq 0$
for any essential edge $e$, and $\omega (e) = 1$ for all other edges.

\begin{theorem}
\label{th:essential weightings2} 
The map $\omega \mapsto x(\Gamma_0,\omega)$ restricts to a
bijection between the set of all essential positive weightings of
$\Gamma_0$ and the set 
of all totally positive $n \times n$ matrices. 
\end{theorem}

The proof of this theorem will use the following notions. 
A minor $\Delta_{I,J}$ is called \emph{solid} if both $I$ and $J$
consist of several consecutive indices;
if furthermore $I\cup J$ contains~$1$,
then $\Delta_{I,J}$ is called \emph{initial} 
(see Figure~\ref{fig:gasca}).   
Each matrix entry is the lower-right corner of exactly one initial minor; 
thus the total number of such minors is~$n^2$. 

\begin{lemma}
\label{lem:initial minors} 
The $n^2$ weights of essential edges in an essential weighting 
$\omega$ of $\Gamma_0$ are related to the 
$n^2$ initial minors of the weight matrix $x=x(\Gamma_0,\omega)$ by 
an invertible monomial transformation.
Thus an essential weighting 
$\omega$ of $\Gamma_0$ is uniquely recovered from~$x$.
\end{lemma}

\proof
The network $\Gamma_0$ has the following easily verified property: 
for any set $I$ of $k$ consecutive indices in $[1,n]$, there is a
unique collection of $k$ vertex-disjoint paths connecting the sources
labeled by~$[1,k]$ (resp.\ by $I$) with the sinks labeled by $I$
(resp.\ by~$[1,k]$). 
These paths are shown by dotted lines in Figure~\ref{fig:maxnetwork},
for $k=2$ and $I=[3,4]$. 
By Lindstr\"om's Lemma, every initial minor~$\Delta$ of 
$x(\Gamma_0,\omega)$ is equal to the product of the weights of
essential edges covered by this family of paths. 
Notice that among these edges, there is always a unique 
\emph{uppermost} essential edge $e(\Delta)$ (indicated by the arrow in
Figure~\ref{fig:maxnetwork}).
Furthermore, the map $\Delta\mapsto e(\Delta)$ is a bijection between
initial minors and essential edges. 
It follows that the weight of each essential edge $e=e(\Delta)$ is equal
to $\Delta$ times a Laurent monomial in some initial minors $\Delta'$
whose associated edges $e(\Delta')$ are located below~$e$. 
%
\endproof

To illustrate Lemma~\ref{lem:initial minors}, consider the special
case~$n=3$. 
The network $\Gamma_0$ is shown in Figure~\ref{fig:planar};
its essential edges have the weights $a, b, \dots, i$. 
The weight matrix $x(\Gamma_0,\omega)$ 
is given in~(\ref{eq:3by3}).
Its initial minors are given by the monomials 
\[
\begin{array}{lll}
\Delta_{1,1}=\underline {d}       & \Delta_{1,2}\ \ \,= d \underline {h}    &
   \Delta_{1,3}\ \ \ \ \,= dh \underline {i}     \\ 
\Delta_{2,1}= \underline {b}  d  & \Delta_{12,12}= d \underline {e}   &
   \Delta_{12,23}\ \ \,= de \underline {g}  h \\
\Delta_{3,1}= \underline {a} bd & \Delta_{23,12}= b\underline {c} de &
   \Delta_{123,123}= de \underline {f}  
\end{array}
\]
where for each minor $\Delta$, the ``leading entry" $\omega(e(\Delta))$ is
underlined. 

To complete the proof of Theorem~\ref{th:essential weightings2}, it
remains to show  
that every totally positive matrix $x$ has the form $x(\Gamma_0,\omega)$
for some essential positive weighting $\omega$. 
By Lemma~\ref{lem:initial minors}, 
such an $\omega$ can be chosen so that $x$ and $x(\Gamma_0,\omega)$ 
will have the same initial minors. 
Thus our claim will follow from the lemma below.

\begin{lemma}
\label{lem:initial minors2}
A square matrix~$x$ is uniquely determined by its initial minors provided
all these minors are nonzero.
\end{lemma}

\proof
Let us show that each matrix entry $x_{ij}$ of~$x$ is
uniquely determined by the initial minors. 
If $i\! =\! 1$ or $j\! =\! 1$, there is nothing to prove,
since $x_{ij}$ is itself an initial minor. 
Assume that $\min (i, j) > 1$.
Let $\Delta$ be the initial minor whose last row is~$i$ and last
column is~$j$, and let $\Delta'$ be the initial minor obtained from
$\Delta$ by deleting this row and this column.  
Then $\Delta= \Delta' x_{ij} + P$, where $P$ is a polynomial 
in the matrix entries $x_{i'j'}$ with $(i',j') \neq (i,j)$ 
and $i' \leq i, \, j' \leq j$. 
Using induction on $i + j$, we can assume that 
each $x_{i'j'}$ that occurs in $P$ is uniquely
determined by the initial minors, 
so the same is true for $x_{ij} = (\Delta - P)/ \Delta'$. 
This completes the proofs of Lemma~\ref{lem:initial minors2} and 
Theorem~\ref{th:essential weightings2}.
\endproof

Theorem~\ref{th:essential weightings2} describes a
\emph{parametrization} of totally positive matrices by 
$n^2$-tuples of positive reals,
providing a partial answer (one of the many possible, as we will see)
to the first problem stated in the introduction. 
The second problem---that of \emph{testing} total
positivity of a matrix---can also be solved using this theorem,
as we will now explain. 

An $n \times n$ matrix has altogether $\binom{2n}{n} - 1$ minors. 
This makes it impractical to test positivity of every single minor. 
It is desirable to find efficient criteria for total positivity
that would only check a small fraction of all minors.

\begin{example}{\rm 
A $2\times 2$ matrix $x = \matbr{a}{b}{c}{d}$ has 
$\binom{4}{2} - 1 = 5$ minors: 
four matrix entries and the determinant $\Delta = ad - bc$.
To test that $x$ is totally positive, it is enough to check positivity
of $a$, $b$, $c$, and~$\Delta$; 
then $d = (\Delta + bc)/a>0$.  
}\end{example}

The following theorem generalizes this example to matrices of arbitrary size;
it is a direct corollary of Theorem~\ref{th:essential weightings2} and 
Lemmas~\ref {lem:initial minors} and~\ref{lem:initial minors2}. 

\begin{theorem}
\label{th:GP-criterion}
A square matrix is totally positive if and only if all its initial minors 
(see Figure~\ref{fig:gasca}) are positive. 
\end{theorem}

This criterion involves $n^2$ minors, 
and it can be shown that this number cannot be lessened.
Theorem~\ref{th:GP-criterion} was proved by M.~Gasca and
J.~M.~Pe\~na~\cite[Theorem~4.1]{gasca-pena}
(for rectangular matrices); 
it also follows from C.~Cryer's results in~\cite{cryer}.
Theorem~\ref{th:GP-criterion} is an enhancement of 
the 1912 criterion by M.~Fekete~\cite{fekete},
who proved that the positivity of all 
solid minors of a matrix implies its total positivity.

\section*{
Theorems of Whitney and Loewner
}
\label{sec:chips}

In this paper, we shall only consider 
\emph{invertible} totally nonnegative $n \times n$ matrices.
Although these matrices have real entries, it is convenient to view
them as elements of the general linear group $G = GL_n (\CC)$. 
We denote by $G_{\geq 0}$ (resp.\ $G_{> 0}$) 
the set of all totally nonnegative (resp.\ totally positive) matrices 
in~$G$.
The structural theory of these matrices begins with the following
basic observation, which is an immediate corollary of 
the Binet-Cauchy formula.

\begin{proposition}
\label{pr:multiplicativity}
Both $G_{\geq 0}$ and $G_{> 0}$ are closed under matrix multiplication.
Furthermore, if $x \in G_{\geq 0}$ and $y \in G_{> 0}$ then 
both $xy$ and $yx$ belong to $G_{> 0}$. 
\end{proposition}

Combining this proposition with the foregoing results, 
we will prove the following theorem of 
A.~Whitney~\cite{whitney}. 

\begin{theorem} 
\label{th:Whitney classical}
{\rm (Whitney's 
Theorem)} 
Every invertible totally nonnegative matrix is the 
limit of a sequence of totally positive matrices. 
\end{theorem}

Thus $G_{\geq 0}$ is the closure of $G_{> 0}$ in~$G$.
(The condition of invertibility in Theorem~\ref{th:Whitney classical} 
can in fact be lifted.) 

\proof 
First let us show that the identity matrix $I$ lies in the closure of
$G_{> 0}\,$. 
By Corollary~\ref{cor:minimal network}, it suffices to show that 
$I = \lim_{N \to \infty} x(\Gamma_0,\omega_N)$ for some sequence of
positive weightings $\omega_N$ of the network $\Gamma_0\,$. 
Note that the map $\omega \mapsto x(\Gamma_0,\omega)$ is continuous,
and choose any sequence 
of positive weightings that converges to the weighting $\omega_0$
defined by $\omega_0 (e) = 1$ 
(resp.~$0$) for all horizontal (resp.\ slanted) edges~$e$. 
Clearly, $x(\Gamma_0,\omega_0) = I$, as desired. 

To complete the proof, write any matrix $x \in G_{\geq 0}$ as  
$x = \lim_{N \to \infty} x \cdot x(\Gamma_0,\omega_N) \,$, 
and notice that all matrices $x \cdot x(\Gamma_0,\omega_N)$ are
totally positive by Proposition~\ref{pr:multiplicativity}. 
\endproof

The following description of the multiplicative monoid $G_{\geq 0}$ 
was first given by C.~Loewner~\cite{loewner}
under the name ``Whitney's Theorem''; it can indeed be deduced 
from~\cite{whitney}. 

\begin{theorem} 
\label{th:Loewner classical}
{\rm (Loewner-Whitney Theorem)} 
\ Any invertible totally nonnegative matrix is 
a product of elementary Jacobi matrices with nonnegative matrix entries.
\end{theorem}

Here ``elementary Jacobi matrix'' is a matrix $x\!\in\! G$ that differs
from~$I$ in a single entry located either on the
main diagonal, or immediately above or below it.

\proof 
We start with an inventory of elementary Jacobi matrices. 
Let $E_{i,j}$ denote the $n \times n$ matrix whose 
$(i,j)$-entry is~1 while all other entries are~0.
For $t\in\CC$ and $i=1,\dots,n-1$, let
\[
x_i(t)=I+tE_{i,i+1}=
\left[\begin{array}{cccccc}
1      & \cdots & 0      & 0      & \cdots & 0      \\
\cdots & \cdots & \cdots & \cdots & \cdots & \cdots \\
0      & \cdots & 1      & t      & \cdots & 0      \\
0      & \cdots & 0      & 1      & \cdots & 0      \\
\cdots & \cdots & \cdots & \cdots & \cdots & \cdots \\
0      & \cdots & 0      & 0      & \cdots & 1 \\
\end{array}\right]
\]
and
\[
x_{\bar i}(t)=I+tE_{i+1,i}=
\left(x_i(t) \right)^T 
\]
(the transpose of $x_i(t)$). 
Also, for $i\!=\!1,\dots, n$ and $t\neq 0$,~let
\[
x_\circi(t)=I+(t-1)E_{i,i} \, ,
\]
the diagonal matrix with $i$th diagonal entry equal to~$t$ and all
other diagonal entries equal to~1. 
Thus elementary Jacobi matrices are precisely the matrices of the form
$x_i(t)$, $x_{\bar i}(t)$, and $x_\circi(t)$.  
An easy check shows that they are totally nonnegative for any~$t > 0$.

\begin{figure}[ht]
\setlength{\unitlength}{1.2pt} 
\begin{center}
\begin{picture}(50,60)(0,-15)
\thicklines

\put(5,0){\line(1,0){20}}
\put(5,20){\line(1,0){20}}
\put(5,40){\line(1,0){20}}
\dark{
\put(5,20){\line(1,1){20}}
}
  \put(5,0){\circle*{2.5}}
  \put(5,20){\circle*{2.5}}
  \put(5,40){\circle*{2.5}}

  \put(25,0){\circle*{2.5}}
  \put(25,20){\circle*{2.5}}
  \put(25,40){\circle*{2.5}}

  \put(7,28){$t$}

  \put(4,-17){$x_i(t)$}

\end{picture}
\begin{picture}(50,60)(0,-15)
\thicklines

\put(5,0){\line(1,0){20}}
\put(5,20){\line(1,0){20}}
\put(5,40){\line(1,0){20}}

  \put(5,0){\circle*{2.5}}
  \put(5,20){\circle*{2.5}}
  \put(5,40){\circle*{2.5}}

  \put(25,0){\circle*{2.5}}
  \put(25,20){\circle*{2.5}}
  \put(25,40){\circle*{2.5}}

\thinlines
\light{
\put(5,40){\line(1,-1){20}}
}
  \put(20,28){$t$}

  \put(4,-17){$x_{\bar i}(t)$}

\end{picture}
\begin{picture}(50,60)(0,-15)
\thicklines

\put(5,0){\line(1,0){20}}
\bulletcolor{
\put(5,20){\line(1,0){20}}
}
\put(5,40){\line(1,0){20}}

  \put(5,0){\circle*{2.5}}
  \put(5,20){\circle*{2.5}}
  \put(5,40){\circle*{2.5}}

  \put(25,0){\circle*{2.5}}
  \put(25,20){\circle*{2.5}}
  \put(25,40){\circle*{2.5}}

  \put(13,24){$t$}

  \put(2,-17){$x_{\circi}(t)$}

\end{picture}
\end{center}
\caption{Elementary ``chips''}
\label{fig:chips}
\end{figure}
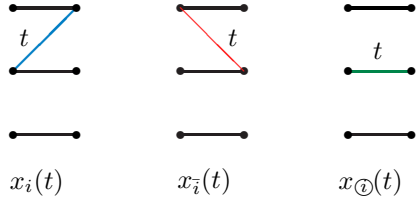

For any word 
$\ii=(i_1, \dots, i_l)$ in the alphabet
\begin{equation}
\label{eq:alphabet}
\A=\{1, \dots, n-1, 
\circled{\,1}\,, \dots, \circled{\,n}\,, 
\overline{1}, \dots, \overline{n-1}\} 
\end{equation}
we define the \emph{product map} $x_\ii : (\CC\setminus\{0\})^l \to G$ by
\begin{equation}
\label{eq:product-map}
x_{\ii}(t_1,\dots,t_l)=x_{i_1}(t_1)\cdots x_{i_l}(t_l) \ .
\end{equation}
(Actually, $x_{\ii}(t_1,\dots,t_l)$ is well defined as long as
the right-hand side of (\ref{eq:product-map}) does not involve  
any factors
of the form $x_{\circi}(0)$.) 
To illustrate, the word 
$\,\ii=\circled{\,1}\, \bar 1 \circled{\,2}\, 1$
gives rise to 
\[
\begin{array}{rcl}
x_\ii(t_1,t_2,t_3,t_4) 
&\!\!\!\!\!=& 
\!\!\!\!\!\matbr{t_1}{0}{0}{1}
\matbr{1}{0}{t_2}{1}
\matbr{1}{0}{0}{t_3}
\matbr{1}{t_4}{0}{1} \\[.2in]
&\!\!\!\!\!=&  \!\!\!\!\!
\matbr{t_1}{t_1t_4}{t_2}{t_2t_4\!+\!t_3} \, .
\end{array}
\]
We will interpret each matrix $x_{\ii}(t_1,\dots,t_l)$ 
as the weight matrix of a planar network.
First note that any elementary Jacobi matrix is the weight matrix of 
a ``chip'' of one of the three kinds shown in Figure~\ref{fig:chips}. 
In each ``chip,'' all edges but one have weight~$1$;
the distinguished edge has weight~$t$. 
Slanted edges connect horizontal levels $i$ and $i+1$, counting from
the bottom; in all examples in Figure~\ref{fig:chips}, $i=2$.

\noindent

The weighted planar network $(\Gamma(\ii), \omega (t_1, \dots, t_l))$
is then constructed by concatenating the ``chips'' 
corresponding to consecutive factors $x_{i_k} (t_k)$,
as shown in Figure~\ref{fig:planar-network}. 
It is easy to see that concatenation of planar networks 
corresponds to multiplying their weight matrices.
We conclude that
the 
product $x_{\ii}(t_1,\dots,t_l)$ of elementary Jacobi matrices equals 
the weight matrix $x(\Gamma (\ii), \omega (t_1,\dots,t_l))$.

\begin{figure}[hb]
\setlength{\unitlength}{1.2pt} 
\begin{center}
\begin{picture}(180,55)(0,-10)
\thicklines

\put(0,40){\line(1,0){40}}
\put(60,40){\line(1,0){120}}
\put(0,20){\line(1,0){160}}
\put(0,0){\line(1,0){100}}
\put(120,0){\line(1,0){60}}

\bulletcolor{
\put(40,40){\line(1,0){20}}
\put(100,0){\line(1,0){20}}
\put(160,20){\line(1,0){20}}
}

\dark{
\put(20,0){\line(1,1){20}}
\put(60,20){\line(1,1){20}}
\put(140,0){\line(1,1){20}}
}

\thinlines
\light{
\put(0,40){\line(1,-1){20}}
\put(80,20){\line(1,-1){20}}
\put(120,40){\line(1,-1){20}}
}

  \put(185,-2){${1}$}
  \put(185,18){${2}$}
  \put(185,38){${3}$}
  \put(-9,-2){${1}$}
  \put(-9,18){${2}$}
  \put(-9,38){${3}$}

  \put(0,0){\circle*{2.5}}
  \put(0,20){\circle*{2.5}}
  \put(0,40){\circle*{2.5}}
  \put(20,0){\circle*{2.5}}
  \put(20,20){\circle*{2.5}}
  \put(20,40){\circle*{2.5}}
  \put(40,0){\circle*{2.5}}
  \put(40,20){\circle*{2.5}}
  \put(40,40){\circle*{2.5}}
  \put(60,0){\circle*{2.5}}
  \put(60,20){\circle*{2.5}}
  \put(60,40){\circle*{2.5}}
  \put(80,0){\circle*{2.5}}
  \put(80,20){\circle*{2.5}}
  \put(80,40){\circle*{2.5}}
  \put(100,0){\circle*{2.5}}
  \put(100,20){\circle*{2.5}}
  \put(100,40){\circle*{2.5}}
  \put(120,0){\circle*{2.5}}
  \put(120,20){\circle*{2.5}}
  \put(120,40){\circle*{2.5}}
  \put(140,0){\circle*{2.5}}
  \put(140,20){\circle*{2.5}}
  \put(140,40){\circle*{2.5}}
  \put(160,0){\circle*{2.5}}
  \put(160,20){\circle*{2.5}}
  \put(160,40){\circle*{2.5}}
  \put(180,0){\circle*{2.5}}
  \put(180,20){\circle*{2.5}}
  \put(180,40){\circle*{2.5}}

\dark{
  \put(35,8){$t_2$}
  \put(75,28){$t_4$}
  \put(155,8){$t_8$}
}

\light{

  \put(15,28){$t_1$}
  \put(80,8){$t_5$}
  \put(135,28){$t_7$}

}

\bulletcolor{
  \put(47,44){$t_3$}
  \put(107,4){$t_6$}
  \put(167,24){$t_9$}
}

\end{picture}
\end{center}
\begin{center}
{$\ii=\overline 2~~ 1\circled{\,3}\,\,2~~\overline 1\,
\circled{\,1}\,\,\overline 2~~1\circled{\,2}$}
\end{center}
\caption{Planar network~$\Gamma(\ii)$}
\label{fig:planar-network}
\end{figure}

In particular, the network $(\Gamma_0, \omega)$ 
appearing in Figure~\ref{fig:maxnetwork} and
Theorem~\ref{th:essential weightings2} 
(more precisely, its equivalent deformation) 
corresponds to some special word $\ii_{\max}$ of length~$n^2$;
instead of defining $\ii_{\max}$ formally, we just write
it for~$n\!=\!4$:  
\[
\ii_{\max}\!=\! (\overline 3,\overline 2, \overline 3, \overline 1,\overline
2,\overline 3, \!\!\circled{\,1},\!\circled{\,2},\!\circled{\,3},\!\circled{\,4},
3,2, 3,1,2,3) \ .
\] 
In view of this,
Theorem~\ref{th:essential weightings2} can be reformulated as follows.  

\begin{theorem}
\label{th:imax parametrization}
The product map $x_{\ii_{\max}}$ restricts to a bijection between 
$n^2$-tuples of positive real numbers and 
totally positive $n \times n$ matrices.
\end{theorem}

We will prove the following refinement of Theorem~\ref{th:Loewner
classical},
which is a reformulation of its original version~\cite{loewner}. 

\begin{theorem}
\label{th:Loewner imax}
Every matrix $x \in G_{\geq 0}$ 
can be written as $x=x_{\ii_{\max}}(t_1, \dots,
t_{n^2})$,  
for some 
$t_1, \dots, t_{n^2}\geq 0$.
\end{theorem}

(Since $x$ is invertible, we must in fact have $t_k > 0$ for 
$\frac{n(n-1)}{2} < k \leq \frac{n(n+1)}{2}$, i.e., for those indices $k$ 
for which the corresponding entry of $\ii_{\max}$ is of the
form~$\circi\,$.)  

\proof
The following key lemma is due to C.~Cryer~\cite{cryer}. 

\begin{lemma}
\label{lem:Gauss}
The leading principal minors $\Delta_{[1,k], [1,k]}$ of a matrix
$x \in G_{\geq 0}$ are positive, for $k\!=\!1,\dots,n$. 
\end{lemma}

\proof
Using induction on $n$, it suffices to show that 
$\Delta_{[1,n-1], [1,n-1]} (x) > 0$.
Let $\Delta^{i,j} (x)$ (resp. $\Delta^{ii',jj'} (x)$) denote 
the minor of~$x$ obtained by deleting the row $i$ and the column $j$
(resp.\ rows $i$ and $i'$, and columns $j$ and $j'$). 
Then, for any $1\leq i<i'\leq n$ and $1\leq j < j'\leq n$, one has  
\begin{equation} 
\label{eq:Dodgson}
\! \Delta^{i'\!,j'\!}(x) \Delta^{i,j}(x) 
-\Delta^{i'\!,j}(x)  \Delta^{i,j'\!} (x)  
=\det (x) \Delta^{ii'\!, jj'\!} (x) 
\end{equation}
as an immediate consequence of Jacobi's formula for minors of the
inverse matrix (see, e.g.,~\cite[Lemma~9.2.10]{brualdi-ryser}). 
The determinantal identity~(\ref{eq:Dodgson}) was proved by
P.~Desnanot as early as in 1819 (see~\cite[pp.~140-142]{muir}); 
it is sometimes called ``Lewis Carroll's identity,'' due to the role
it plays in C.~L.~Dodgson's condensation
method \cite[pp.~170--180]{dodgson}. 

Now suppose that $\Delta^{n,n} (x) = 0$ for some $x \in G_{\geq 0}\,$. 
Since $x$ is invertible, we have $\Delta^{i,n}(x)>0$ and 
$\Delta^{n,j}(x)>0$ for some indices $i, j < n$. 
Using (\ref{eq:Dodgson}) with $i' = j' = n$, we arrive at a desired
contradiction by 
\[
0>- \Delta^{n,j}(x) \Delta^{i,n}(x) =  
\det(x) \Delta^{in,jn}(x) \geq 0\,. \endproofmath 
\]

We are now ready to complete the proof of Theorem~\ref{th:Loewner imax}.
Any matrix~$x \in G_{\geq 0}$ is by Theorem~\ref{th:Whitney
classical} a limit of totally positive matrices~$x_N\,$,
each of which can by Theorem~\ref{th:imax parametrization}
be factored as 
$x_N=x_{\ii_{\max}}(t_1^{(N)},\dots,t_{n^2}^{(N)})$ with all
$t_k^{(N)}$ positive. 
It suffices to show that the sequence $s_N = \sum_{k=1}^{n^2}
t_k^{(N)}$ converges; 
then the standard compactness argument will imply that
the sequence of vectors $(t_1^{(N)},\dots,t_{n^2}^{(N)})$ contains a
converging subsequence, whose limit $(t_1,\dots,t_{n^2})$ will provide the
desired factorization $x=x_{\ii_{\max}}(t_1,\dots,t_{n^2})$.
To see that $(s_N)$ converges, we use the explicit formula
\[
\begin{array}{l}
s_N = \displaystyle \sum_{i = 1}^n \frac{\Delta_{[1,i], [1,i]} (x_N)}
{\Delta_{[1,i-1], [1,i-1]} (x_N)}  \\[.2in]
\ \ + 
\displaystyle \sum_{i = 1}^{n-1} \frac
{\Delta_{[1,i-1] \cup \{i+1\}, [1,i]} (x_N) + \Delta_{[1,i], [1,i-1]
    \cup \{i+1\}} (x_N)}  
{\Delta_{[1,i], [1,i]} (x_N)}
\end{array}
\]
(to prove this, compute the minors on the right with the help of 
Lindstr\"om's Lemma and simplify). 
Thus $s_N$ is expressed as a Laurent polynomial in the minors of $x_N$
whose denominators only involve leading principal minors
$\Delta_{[1,k], [1,k]}$.  
By Lemma~\ref{lem:Gauss}, as $x_N$ converges to $x$, 
this Laurent polynomial converges to its value at $x$. 
This completes the proofs of Theorems~\ref{th:Loewner classical}
and \ref{th:Loewner imax}.
\endproof 

\section*{Double wiring diagrams and \\
total positivity criteria}
\label{sec:double}

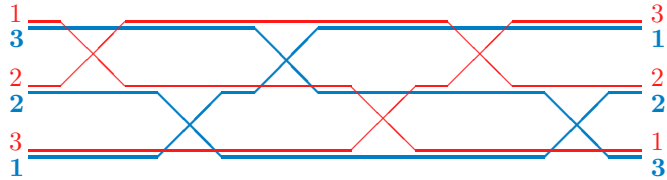
\begin{figure}[ht]
\setlength{\unitlength}{1.22pt} 
\begin{center}
\begin{picture}(180,45)(6,0)
\thicklines
\dark{

  \put(0,0){\line(1,0){40}}
  \put(60,0){\line(1,0){100}}
  \put(180,0){\line(1,0){10}}
  \put(0,20){\line(1,0){40}}
  \put(60,20){\line(1,0){10}}
  \put(90,20){\line(1,0){70}}
  \put(180,20){\line(1,0){10}}
  \put(0,40){\line(1,0){70}}
  \put(90,40){\line(1,0){100}}

  \put(40,0){\line(1,1){20}}
  \put(70,20){\line(1,1){20}}
  \put(160,0){\line(1,1){20}}

  \put(40,20){\line(1,-1){20}}
  \put(70,40){\line(1,-1){20}}
  \put(160,20){\line(1,-1){20}}
  \put(193,-6){$\mathbf{3}$}
  \put(193,14){$\mathbf{2}$}
  \put(193,34){$\mathbf{1}$}

  \put(-6,-6){$\mathbf{1}$}
  \put(-6,14){$\mathbf{2}$}
  \put(-6,34){$\mathbf{3}$}
}

\light{
\thinlines

  \put(0,2){\line(1,0){100}}
  \put(120,2){\line(1,0){70}}
  \put(0,22){\line(1,0){10}}
  \put(30,22){\line(1,0){70}}
  \put(120,22){\line(1,0){10}}
  \put(150,22){\line(1,0){40}}
  \put(0,42){\line(1,0){10}}
  \put(30,42){\line(1,0){100}}
  \put(150,42){\line(1,0){40}}

  \put(10,22){\line(1,1){20}}
  \put(100,2){\line(1,1){20}}
  \put(130,22){\line(1,1){20}}

  \put(10,42){\line(1,-1){20}}
  \put(100,22){\line(1,-1){20}}
  \put(130,42){\line(1,-1){20}}

  \put(193,2){${1}$}
  \put(193,22){${2}$}
  \put(193,42){${3}$}

  \put(-6,2){${3}$}
  \put(-6,22){${2}$}
  \put(-6,42){${1}$}
}
\end{picture}
\end{center}
\caption{Double wiring diagram}
\label{fig:double-wiring}
\end{figure}

We will now give another proof of Theorem~\ref{th:GP-criterion},  
which will include it into a family of ``optimal" 
total positivity criteria that
correspond to combinatorial objects called
\emph{double wiring diagrams}. 
This notion is best explained by an
example, such as the one given in Figure~\ref{fig:double-wiring}.
A double wiring diagram consists of two families of $n$
piecewise-straight 
lines (each family colored with one of the two colors), 
the crucial requirement being that each
pair of 
lines of like color intersect exactly once. 

\begin{figure}[ht]
\setlength{\unitlength}{1.22pt} 
\begin{center}
\begin{picture}(180,60)(6,-10)
\thicklines
\dark{

  \put(0,0){\line(1,0){40}}
  \put(60,0){\line(1,0){100}}
  \put(180,0){\line(1,0){10}}
  \put(0,20){\line(1,0){40}}
  \put(60,20){\line(1,0){10}}
  \put(90,20){\line(1,0){70}}
  \put(180,20){\line(1,0){10}}
  \put(0,40){\line(1,0){70}}
  \put(90,40){\line(1,0){100}}

  \put(40,0){\line(1,1){20}}
  \put(70,20){\line(1,1){20}}
  \put(160,0){\line(1,1){20}}

  \put(40,20){\line(1,-1){20}}
  \put(70,40){\line(1,-1){20}}
  \put(160,20){\line(1,-1){20}}

  \put(193,-6){$\mathbf{3}$}
  \put(193,14){$\mathbf{2}$}
  \put(193,34){$\mathbf{1}$}

  \put(-6,-6){$\mathbf{1}$}
  \put(-6,14){$\mathbf{2}$}
  \put(-6,34){$\mathbf{3}$}
}

\light{
\thinlines

  \put(0,2){\line(1,0){100}}
  \put(120,2){\line(1,0){70}}
  \put(0,22){\line(1,0){10}}
  \put(30,22){\line(1,0){70}}
  \put(120,22){\line(1,0){10}}
  \put(150,22){\line(1,0){40}}
  \put(0,42){\line(1,0){10}}
  \put(30,42){\line(1,0){100}}
  \put(150,42){\line(1,0){40}}

  \put(10,22){\line(1,1){20}}
  \put(100,2){\line(1,1){20}}
  \put(130,22){\line(1,1){20}}

  \put(10,42){\line(1,-1){20}}
  \put(100,22){\line(1,-1){20}}
  \put(130,42){\line(1,-1){20}}

  \put(193,2){${1}$}
  \put(193,22){${2}$}
  \put(193,42){${3}$}

  \put(-6,2){${3}$}
  \put(-6,22){${2}$}
  \put(-6,42){${1}$}

  \put(94,-10){$_{\light{\emptyset},\dark{\emptyset\hspace{-.057in}
\emptyset\hspace{-.057in}\emptyset}}$}

  \put(14,10){$_{\light{3},\mathbf{\dark{1}}}$}
  \put(76,10){$_{\light{3},\mathbf{\dark{2}}}$}
  \put(136,10){$_{\light{1},\mathbf{\dark{2}}}$}
  \put(181,10){$_{\light{1},\mathbf{\dark{3}}}$}

  \put(-1,30){$_{\light{23},\mathbf{\dark{12}}}$}
  \put(42,30){$_{\light{13},\mathbf{\dark{12}}}$}
  \put(102,30){$_{\light{13},\mathbf{\dark{23}}}$}
  \put(162,30){$_{\light{12},\mathbf{\dark{23}}}$}

  \put(84,50){$_{\light{123},\mathbf{\dark{123}}}$}
}
\end{picture}
\end{center}
\caption{Chamber minors}
\label{fig:chamber-sets}
\end{figure}
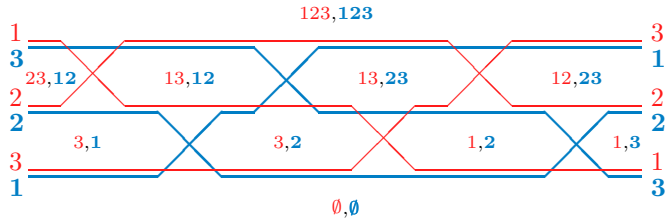

The lines in a double wiring diagram are numbered separately within
each color. 
We then assign to every \emph{chamber} of a diagram a pair of subsets
of the set $[1,n]=\{1,\dots,n\}$: 
each subset indicates which lines of the corresponding
color pass \emph{below} that chamber; 
see Figure~\ref{fig:chamber-sets}. 

Thus every chamber is naturally associated with a minor $\Delta_{I,J}$ 
of an $n\times n$ matrix~$x=(x_{ij})$ (we call it a \emph{chamber minor})
that occupies the rows and columns 
specified by the sets $I$ and $J$ written into that chamber. 
In our running example, there are 9 chamber minors 
(the total number is always~$n^2$),
namely $x_{31}\,$, $x_{32}\,$, $x_{12}\,$, $x_{13}\,$, 
$\Delta_{23,12}\,$, $\Delta_{13,12}\,$, $\Delta_{13,23}\,$,
$\Delta_{12,23}\,$, and $\Delta_{123,123}=\det(x)$. 

\begin{theorem}
\label{th:double diagram TP-criterion}
\cite{FZ}
Every double wiring diagram gives rise to the following 
criterion: an $n\times n$ matrix is totally positive if and
only if all its $n^2$ chamber minors are positive. 
\end{theorem}

\begin{figure}[ht]
\setlength{\unitlength}{1.22pt} 
\begin{center}
\begin{picture}(180,55)(6,0)
\thicklines
\dark{

  \put(0,0){\line(1,0){100}}
  \put(120,0){\line(1,0){40}}
  \put(180,0){\line(1,0){10}}
  \put(0,20){\line(1,0){100}}
  \put(120,20){\line(1,0){10}}
  \put(150,20){\line(1,0){10}}
  \put(180,20){\line(1,0){10}}
  \put(0,40){\line(1,0){130}}
  \put(150,40){\line(1,0){40}}

  \put(100,0){\line(1,1){20}}
  \put(130,20){\line(1,1){20}}
  \put(160,0){\line(1,1){20}}

  \put(100,20){\line(1,-1){20}}
  \put(130,40){\line(1,-1){20}}
  \put(160,20){\line(1,-1){20}}

  \put(193,-6){$\mathbf{3}$}
  \put(193,14){$\mathbf{2}$}
  \put(193,34){$\mathbf{1}$}

  \put(-6,-6){$\mathbf{1}$}
  \put(-6,14){$\mathbf{2}$}
  \put(-6,34){$\mathbf{3}$}
}

\light{
\thinlines

  \put(0,2){\line(1,0){10}}
  \put(30,2){\line(1,0){40}}
  \put(90,2){\line(1,0){100}}
  \put(0,22){\line(1,0){10}}
  \put(30,22){\line(1,0){10}}
  \put(60,22){\line(1,0){10}}
  \put(90,22){\line(1,0){100}}
  \put(0,42){\line(1,0){40}}
  \put(60,42){\line(1,0){130}}

  \put(10,2){\line(1,1){20}}
  \put(40,22){\line(1,1){20}}
  \put(70,2){\line(1,1){20}}

  \put(10,22){\line(1,-1){20}}
  \put(40,42){\line(1,-1){20}}
  \put(70,22){\line(1,-1){20}}

  \put(193,2){${1}$}
  \put(193,22){${2}$}
  \put(193,42){${3}$}

  \put(-6,2){${3}$}
  \put(-6,22){${2}$}
  \put(-6,42){${1}$}

  \put(2,10){$_{\light{3},\mathbf{\dark{1}}}$}
  \put(46,10){$_{\light{2},\mathbf{\dark{1}}}$}
  \put(91,10){$_{\light{1},\mathbf{\dark{1}}}$}
  \put(136,10){$_{\light{1},\mathbf{\dark{2}}}$}
  \put(181,10){$_{\light{1},\mathbf{\dark{3}}}$}

  \put(12,30){$_{\light{23},\mathbf{\dark{12}}}$}
  \put(87,30){$_{\light{12},\mathbf{\dark{12}}}$}
  \put(162,30){$_{\light{12},\mathbf{\dark{23}}}$}

  \put(84,50){$_{\light{123},\mathbf{\dark{123}}}$}
}
\end{picture}
\end{center}
\caption{Lexicographically minimal diagram}
\label{fig:lexmindiagram}
\end{figure}
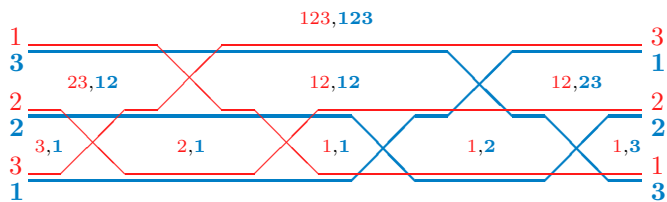

The criterion in Theorem~\ref{th:GP-criterion} 
is a special case of Theorem~\ref{th:double diagram TP-criterion}, 
and arises from the ``lexicographically minimal'' double wiring
diagram, shown in Figure~\ref{fig:lexmindiagram} for~$n=3$. 

\proof
We will actually prove the following statement that implies
Theorem~\ref{th:double diagram TP-criterion}.

\begin{theorem}
\label{th:subtraction-free}
\cite{FZ}
Every minor of a generic square matrix can be written as a 
rational expression in the chamber minors of a given double wiring
diagram, 
and moreover this rational expression is \emph{subtraction-free,}
i.e., all coefficients in the numerator and denominator are positive.
\end{theorem}


Two double wiring diagrams are called \emph{isotopic}
if they have the same collections of chamber minors. 
The terminology suggests what is really going on here:
two isotopic diagrams have the same ``topology.''  From now on, 
we will treat such diagrams as indistinguishable from each other. 

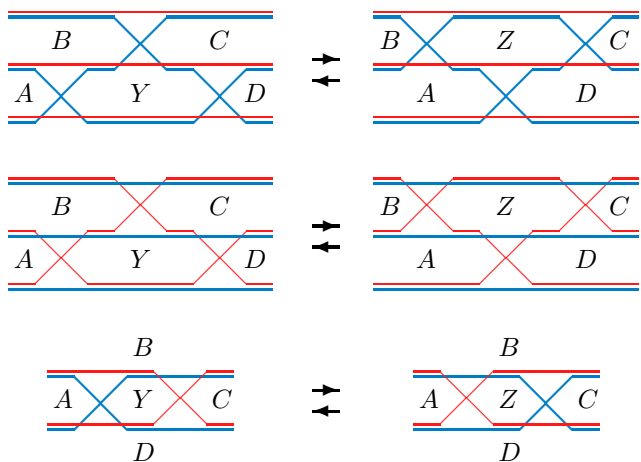
\begin{figure}[ht]
\setlength{\unitlength}{1pt} 
\begin{center}
\begin{picture}(120,45)(10,0)
\thicklines
\dark{

  \put(0,0){\line(1,0){10}}
  \put(30,0){\line(1,0){40}}
  \put(50,0){\line(1,0){10}}
  \put(90,0){\line(1,0){10}}
  \put(0,20){\line(1,0){10}}
  \put(30,20){\line(1,0){10}}
  \put(60,20){\line(1,0){10}}
  \put(90,20){\line(1,0){10}}
  \put(0,40){\line(1,0){40}}
  \put(60,40){\line(1,0){40}}

  \put(10,0){\line(1,1){20}}
  \put(40,20){\line(1,1){20}}
  \put(70,0){\line(1,1){20}}

  \put(10,20){\line(1,-1){20}}
  \put(40,40){\line(1,-1){20}}
  \put(70,20){\line(1,-1){20}}
}
  \put(115,24){\vector(1,0){10}}
  \put(125,16){\vector(-1,0){10}}

\thinlines
\light{
  \put(0,2){\line(1,0){100}}
  \put(0,22){\line(1,0){100}}
  \put(0,42){\line(1,0){100}}
}


  \put(2,8){$A$}
  \put(46,8){$Y$}
  \put(89,8){$D$}

  \put(16,28){$B$}
  \put(76,28){$C$}

\end{picture}
\begin{picture}(100,40)(-5,0)
\thicklines
\dark{

  \put(0,40){\line(1,0){10}}
  \put(30,40){\line(1,0){40}}
  \put(50,40){\line(1,0){10}}
  \put(90,40){\line(1,0){10}}
  \put(0,20){\line(1,0){10}}
  \put(30,20){\line(1,0){10}}
  \put(60,20){\line(1,0){10}}
  \put(90,20){\line(1,0){10}}
  \put(0,0){\line(1,0){40}}
  \put(60,0){\line(1,0){40}}

  \put(10,20){\line(1,1){20}}
  \put(40,0){\line(1,1){20}}
  \put(70,20){\line(1,1){20}}

  \put(10,40){\line(1,-1){20}}
  \put(40,20){\line(1,-1){20}}
  \put(70,40){\line(1,-1){20}}
}

\thinlines
\light{
  \put(0,2){\line(1,0){100}}
  \put(0,22){\line(1,0){100}}
  \put(0,42){\line(1,0){100}}
}
  \put(2,28){$B$}
  \put(46,28){$Z$}
  \put(89,28){$C$}

  \put(16,8){$A$}
  \put(76,8){$D$}

\end{picture}
\end{center}
\begin{center}
\begin{picture}(120,60)(10,0)
\thinlines
\light{

  \put(0,2){\line(1,0){10}}
  \put(30,2){\line(1,0){40}}
  \put(50,2){\line(1,0){10}}
  \put(90,2){\line(1,0){10}}
  \put(0,22){\line(1,0){10}}
  \put(30,22){\line(1,0){10}}
  \put(60,22){\line(1,0){10}}
  \put(90,22){\line(1,0){10}}
  \put(0,42){\line(1,0){40}}
  \put(60,42){\line(1,0){40}}

  \put(10,2){\line(1,1){20}}
  \put(40,22){\line(1,1){20}}
  \put(70,2){\line(1,1){20}}

  \put(10,22){\line(1,-1){20}}
  \put(40,42){\line(1,-1){20}}
  \put(70,22){\line(1,-1){20}}
}
  \put(2,8){$A$}
  \put(46,8){$Y$}
  \put(89,8){$D$}

  \put(16,28){$B$}
  \put(76,28){$C$}

\thicklines
  \put(115,24){\vector(1,0){10}}
  \put(125,16){\vector(-1,0){10}}
\dark{
  \put(0,0){\line(1,0){100}}
  \put(0,20){\line(1,0){100}}
  \put(0,40){\line(1,0){100}}
}

\end{picture}
\begin{picture}(100,40)(-5,0)
\thinlines
\light{

  \put(0,42){\line(1,0){10}}
  \put(30,42){\line(1,0){40}}
  \put(50,42){\line(1,0){10}}
  \put(90,42){\line(1,0){10}}
  \put(0,22){\line(1,0){10}}
  \put(30,22){\line(1,0){10}}
  \put(60,22){\line(1,0){10}}
  \put(90,22){\line(1,0){10}}
  \put(0,2){\line(1,0){40}}
  \put(60,2){\line(1,0){40}}

  \put(10,22){\line(1,1){20}}
  \put(40,2){\line(1,1){20}}
  \put(70,22){\line(1,1){20}}

  \put(10,42){\line(1,-1){20}}
  \put(40,22){\line(1,-1){20}}
  \put(70,42){\line(1,-1){20}}
}

\thicklines
\dark{
  \put(0,0){\line(1,0){100}}
  \put(0,20){\line(1,0){100}}
  \put(0,40){\line(1,0){100}}
}
  \put(2,28){$B$}
  \put(46,28){$Z$}
  \put(89,28){$C$}

  \put(16,8){$A$}
  \put(76,8){$D$}

\end{picture}
\end{center}
\begin{center}
\begin{picture}(100,60)(5,-10)
\thicklines
\dark{

  \put(0,0){\line(1,0){10}}
  \put(30,0){\line(1,0){40}}

  \put(0,20){\line(1,0){10}}
  \put(30,20){\line(1,0){40}}

  \put(10,0){\line(1,1){20}}
  \put(10,20){\line(1,-1){20}}
}

\thinlines
\light{

  \put(0,2){\line(1,0){40}}
  \put(60,2){\line(1,0){10}}

  \put(0,22){\line(1,0){40}}
  \put(60,22){\line(1,0){10}}

  \put(40,2){\line(1,1){20}}
  \put(40,22){\line(1,-1){20}}
}
  \put(2,8){$A$}
  \put(32,8){$Y$}
  \put(62,8){$C$}

  \put(32,28){$B$}
  \put(32,-12){$D$}

\thicklines
  \put(100,15){\vector(1,0){10}}
  \put(110,7){\vector(-1,0){10}}

\end{picture}
\begin{picture}(100,60)(-30,-10)
\thicklines
\dark{
  \put(0,0){\line(1,0){40}}
  \put(60,0){\line(1,0){10}}

  \put(0,20){\line(1,0){40}}
  \put(60,20){\line(1,0){10}}

  \put(40,0){\line(1,1){20}}
  \put(40,20){\line(1,-1){20}}
}

\thinlines
\light{

  \put(0,2){\line(1,0){10}}
  \put(30,2){\line(1,0){40}}

  \put(0,22){\line(1,0){10}}
  \put(30,22){\line(1,0){40}}

  \put(10,2){\line(1,1){20}}
  \put(10,22){\line(1,-1){20}}
}
  \put(2,8){$A$}
  \put(32,8){$Z$}
  \put(62,8){$C$}

  \put(32,28){$B$}
  \put(32,-12){$D$}

\end{picture}
\end{center}
\caption{Local ``moves''}
\label{fig:moves}
\end{figure}

We will deduce Theorem~\ref{th:subtraction-free} from the following fact:
any two double wiring diagrams can be
transformed into each other by a sequence of local ``moves'' of
three different kinds, shown in Figure~\ref{fig:moves}. 
(This is a direct corollary of a theorem of G.~Ringel~\cite{ringel}. 
It can also be derived from the Tits theorem on reduced words in the
symmetric group;
cf.\ (\ref{eq:iji=jij})--(\ref{eq:ab=ba}) below.) 

Notice that each local move exchanges a single chamber minor~$Y$ with another
chamber minor~$Z$, and keeps all other chamber minors in place. 

\begin{lemma}
\label{lem:3-term}
Whenever two double wiring diagrams differ by a single local move of
one of the three types shown in Figure~\ref{fig:moves},
the chamber minors 
appearing there satisfy the identity $AC+BD=YZ$. 
\end{lemma}

The three-term determinantal identities of Lemma~\ref{lem:3-term} 
are well known, although not in this disguised form. 
The last of these identities is nothing but the identity
(\ref{eq:Dodgson}),  
applied to various submatrices of an $n \times n$ matrix. 
The identities corresponding to the top two ``moves'' in
Figure~\ref{fig:moves} are special instances of the classical
Grassmann-Pl\"ucker relations (see,
e.g.,~\cite[(15.53)]{fulton-harris}), and were obtained by P.~Desnanot 
alongside (\ref{eq:Dodgson}) in the same 1819 publication 
that we mentioned before. 

Theorem~\ref{th:subtraction-free} 
is now proved as follows. 
We first note that any minor appears as chamber minor in  
some double wiring diagram. 
Therefore, it suffices to show that the
chamber minors of one diagram can be writen as subtract\-ion-free
rational expressions in the chamber minors of any other diagram. 
This is a direct corollary of Lemma~\ref{lem:3-term}
combined with the fact that any two diagrams 
are related by a sequence of local moves:
indeed, each local move replaces $Y$ by $(AC+BD)/Z$, or $Z$ by
$(AC+BD)/Y$.  
\endproof

\begin{figure}[t]
\setlength{\unitlength}{1.45pt} 
\begin{center}
\begin{picture}(200,350)(-100,-175)


\put(0,-165){\line(2,3){110}}
\put(0,-165){\line(-2,3){110}}
\put(0,-165){\line(0,1){55}}

\put(0,-110){\line(0,1){30}}
\put(0,-110){\line(2,1){40}}
\put(0,-110){\line(-2,1){40}}

\put(40,-90){\line(2,3){40}}
\put(40,-90){\line(0,1){30}}
\put(40,-90){\line(-2,1){40}}

\put(0,-80){\line(-2,1){40}}
\put(0,-80){\line(2,1){40}}

\put(-40,-90){\line(2,1){40}}
\put(-40,-90){\line(0,1){30}}
\put(-40,-90){\line(-2,3){40}}
 
\put(0,-70){\line(0,1){30}}

\put(40,-60){\line(-2,1){40}}
\put(-40,-60){\line(2,1){40}}
\put(40,-60){\line(0,1){30}}
\put(-40,-60){\line(0,1){30}}

\put(0,-40){\line(0,1){20}}
 
\put(0,-20){\line(1,1){20}}
\put(0,-20){\line(-1,1){20}}
\put(40,-30){\line(2,3){20}}
\put(40,-30){\line(-2,3){20}}
\put(80,-30){\line(1,1){30}}
\put(80,-30){\line(-2,3){20}}
\put(-40,-30){\line(2,3){20}}
\put(-40,-30){\line(-2,3){20}}
\put(-80,-30){\line(2,3){20}}
\put(-80,-30){\line(-1,1){30}}


\put(0,165){\line(2,-3){110}}
\put(0,165){\line(-2,-3){110}}
\put(0,165){\line(0,-1){55}}

\put(0,110){\line(0,-1){30}}
\put(0,110){\line(2,-1){40}}
\put(0,110){\line(-2,-1){40}}

\put(40,90){\line(2,-3){40}}
\put(40,90){\line(0,-1){30}}
\put(40,90){\line(-2,-1){40}}

\put(0,80){\line(-2,-1){40}}
\put(0,80){\line(2,-1){40}}

\put(-40,90){\line(2,-1){40}}
\put(-40,90){\line(0,-1){30}}
\put(-40,90){\line(-2,-3){40}}
 
\put(0,70){\line(0,-1){30}}

\put(40,60){\line(-2,-1){40}}
\put(-40,60){\line(2,-1){40}}
\put(40,60){\line(0,-1){30}}
\put(-40,60){\line(0,-1){30}}

\put(0,40){\line(0,-1){20}}
 
\put(0,20){\line(1,-1){20}}
\put(0,20){\line(-1,-1){20}}
\put(40,30){\line(2,-3){20}}
\put(40,30){\line(-2,-3){20}}
\put(80,30){\line(1,-1){30}}
\put(80,30){\line(-2,-3){20}}
\put(-40,30){\line(2,-3){20}}
\put(-40,30){\line(-2,-3){20}}
\put(-80,30){\line(2,-3){20}}
\put(-80,30){\line(-1,-1){30}}

{\small 
\put(3,-22){$abcG$} 
\put(-17,-1){$acFG$} 
\put(-37,-32){$ceFG$} 
\put(-34,-61){$cdeG$} 
\put(2,-40){$bcdG$} 
\put(42,-58){$bdfG$} 
\put(12,-31){$bfEG$} 
\put(23,-1){$abEG$} 
\put(-9,-88){$defG$} 
\put(2,-69){$bcdA$} 
\put(-34,-92){$cdeA$} 
\put(2,-114){$defA$} 
\put(15,-93){$bdf\!A$} 
\put(2,-167){$efgA$} 
\put(-107,-1){$egAB$} 
\put(-77,-32){$ceAB$} 
\put(-57,-1){$ceBF$} 
\put(64,1){$bfCE$} 
\put(53,-33){$bfAC$} 
\put(80,-6){$fgAC$} 
      
\put(2,20){$aEFG$}
 
\put(-77,29){$egBF$} 
\put(-37,29){$acBF$}
\put(14,29){$abCE$} 
\put(51,28){$fgCE$} 
 
\put(3,37){$aDEF$} 

\put(-60,49){$aBDF$} 
\put(28,50){$aCDE$}
 
\put(-13.5,62){$gD\,EF$} 

\put(-13,82){$aB\,CD$}
 
\put(-34,89){$gBDF$} 
\put(8.5,90.5){$gCDE$} 

\put(2,112){$gBCD$} 
\put(2,165){$gABC$} 
      
 
\put(0,-20){\circle*{3}} 
\put(-20,0){\circle*{3}} 
\put(-40,-30){\circle*{3}} 
\put(-40,-60){\circle*{3}} 
\put(0,-40){\circle*{3}} 
\put(40,-60){\circle*{3}} 
\put(40,-30){\circle*{3}} 
\put(20,0){\circle*{3}} 
\put(0,-80){\circle*{3}} 
\put(0,-70){\circle*{3}} 
\put(-40,-90){\circle*{3}} 
\put(0,-110){\circle*{3}} 
\put(40,-90){\circle*{3}} 
\put(0,-165){\circle*{3}} 
\put(-110,0){\circle*{3}} 
\put(-80,-30){\circle*{3}} 
\put(-60,0){\circle*{3}} 
\put(60,0){\circle*{3}} 
\put(80,-30){\circle*{3}} 
\put(110,0){\circle*{3}} 
      
\put(0,20){\circle*{3}} 
 
\put(-80,30){\circle*{3}} 
\put(-40,30){\circle*{3}} 
\put(40,30){\circle*{3}} 
\put(80,30){\circle*{3}} 
 
\put(0,40){\circle*{3}} 

\put(-40,60){\circle*{3}} 
\put(40,60){\circle*{3}} 
 
\put(0,70){\circle*{3}} 

\put(0,80){\circle*{3}} 
 
\put(-40,90) {\circle*{3}} 
\put(40,90){\circle*{3}} 

\put(0,110){\circle*{3}} 
\put(0,165){\circle*{3}}


{\small

\put(-123,-128){
$\begin{array}{ll}
a=x_{11} \\[.05in]
b=x_{12} \\[.05in]
c=x_{21} \\[.05in]
d=x_{22} \\[.05in]
e=x_{23} \\[.05in]
f=x_{32} \\[.05in]
g=x_{33} 
\end{array}
$}
 
\put(60,-128){
$\begin{array}{ll}
A=\Delta_{23,23}\\[.05in]
B=\Delta_{23,13}\\[.05in] 
C=\Delta_{13,23}\\[.05in] 
D=\Delta_{13,13}\\[.05in] 
E=\Delta_{13,12}\\[.05in] 
F=\Delta_{12,13}\\[.05in] 
G=\Delta_{12,12} 
\end{array}
$}

} 

} 

\end{picture}
\end{center}
\caption{
Total positivity criteria for $GL_3$}
\label{fig:schemes}
\end{figure}
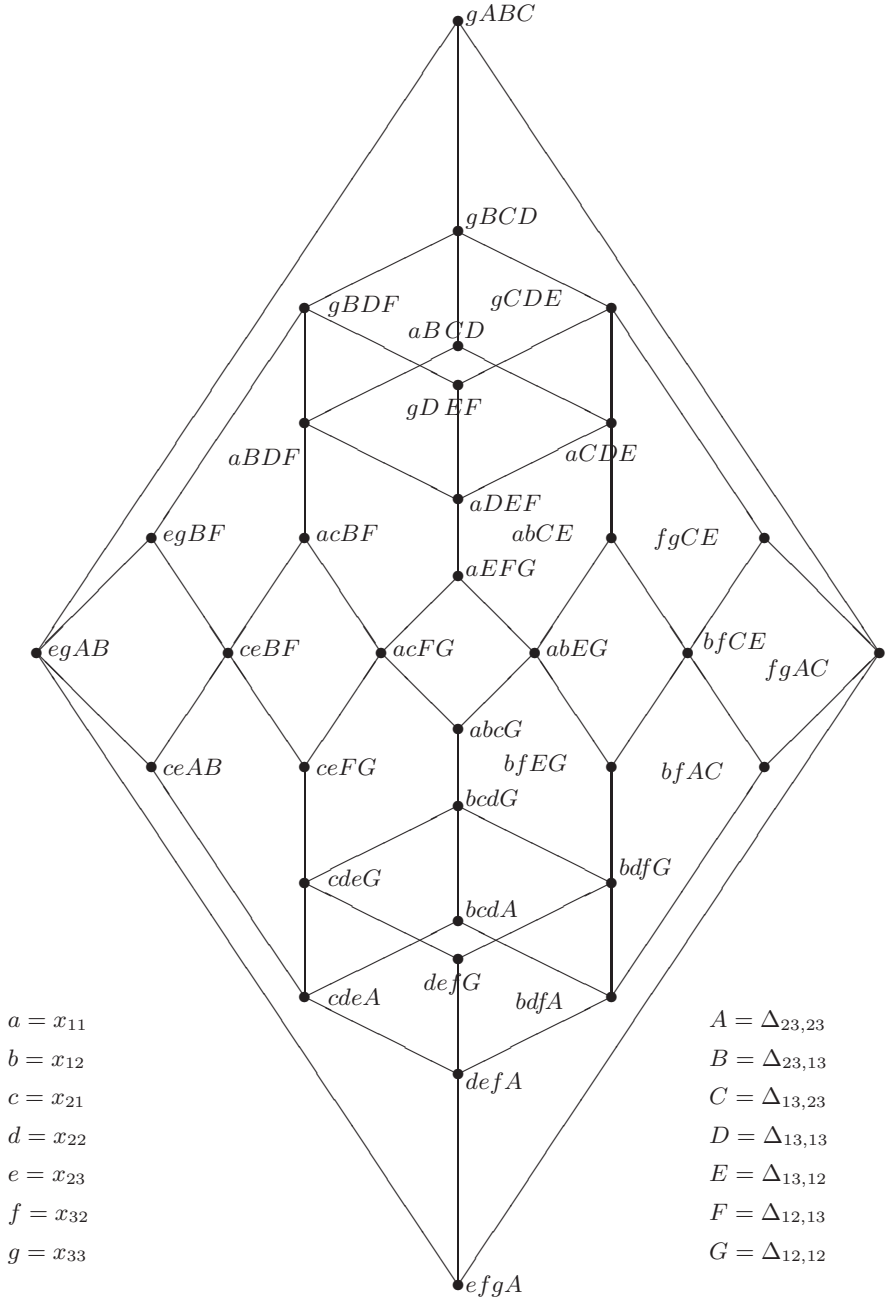

Implicit in the above proof is an important combinatorial structure
lying behind Theorems~\ref{th:double diagram
  TP-criterion} and~\ref{th:subtraction-free}: the graph   
$\Phi_n$ whose vertices are the (isotopy classes of) double wiring
diagrams, and whose edges correspond to local moves. 
The study of $\Phi_n$ is an interesting problem in itself.
The first non\-trivial example is the graph $\Phi_3$ shown in
Figure~\ref{fig:schemes}. 
It has $34$ vertices, corresponding to 
$34$ different total positivity criteria.
Each of these criteria tests 9 minors of a $3\! \times\! 3$ matrix.
Five of these minors, viz.\ 
$x_{31} \,$, 
$x_{13}\,$,
$\Delta_{23,12} \,$, 
$\Delta_{12,23} \,$, 
and $\det(x)$, correspond to the ``unbounded'' chambers 
that lie on the periphery of every double wiring diagram; 
they are common to all 34 criteria. 
The other four minors correspond to the bounded chambers,
and depend on the choice of a diagram. 
For example, the criterion derived from
Figure~\ref{fig:chamber-sets} involves ``bounded''
chamber minors 
$\Delta_{3,2}\,$, $\Delta_{1,2}\,$, $\Delta_{13,12}\,$,  and
$\Delta_{13,23}\,$. 
In Figure~\ref{fig:schemes}, each vertex of $\Phi_3$ is 
labeled by the quadruple of ``bounded'' minors
that appear in the corresponding total positivity criterion.

We suggest the following refinement of 
Theorem~\ref{th:subtraction-free}. 

\begin{conjecture}
\label{conj:Laurent} 
Every minor of a generic square matrix 
can be written as 
a Laurent polynomial
with nonnegative integer coefficients in the chamber minors of an
arbitrary double wiring diagram. 
\end{conjecture}

Perhaps more important than proving this conjecture would be to give
explicit combinatorial expressions for the Laurent polynomials in question.  
We note a case where the conjecture is true, and the desired expressions
can be given: the ``lexicographically minimal'' double wiring
diagram whose chamber minors are the initial minors.
Indeed, a generic matrix $x$ can be uniquely written as the product 
$x_{\ii_{\max}}(t_1, \dots, t_{n^2})$ of elementary Jacobi matrices 
(cf.\ Theorem~\ref{th:imax parametrization});
then each minor of $x$ can be written as a polynomial 
in the $t_k$ with nonnegative integer coefficients (with
the help of Lindstr\"om's Lemma), while each $t_k$ is a Laurent
monomial in the initial minors of~$x$, 
by Lemma~\ref{lem:initial minors}. 

It is proved in~\cite[Theorem~1.13]{FZ} that every minor can be
written as a Laurent polynomial with integer (possibly negative)
coefficients in the chamber minors of a given diagram.
Note, however, that this result, combined with
Theorem~\ref{th:subtraction-free}, does not 
imply Conjecture~\ref{conj:Laurent}, 
because there do exist subtraction-free rational 
expressions that are Laurent polynomials although \emph{not} 
with nonnegative coefficients 
(e.g., think of $(p^3+q^3)/(p+q)=p^2-pq+q^2$).

The following special case of
Conjecture \ref{conj:Laurent} can be derived
from~\cite[Thm.~3.7.4]{BFZ}.

\begin{theorem}
Conjecture~\ref{conj:Laurent} holds 
for all 
wiring diagrams
in which all
intersections of one color precede the intersections of another color. 
\end{theorem}

We do not know an elementary proof of this result;
the proof in~\cite{BFZ} depends on the theory of canonical bases
for quantum general linear groups. 

\subsection*{Digression: Somos sequences}

The three-term relation $AC+BD=YZ$ is surrounded by some magic that
eludes our comprehension. 
We cannot resist mentioning the related problem involving the
\emph{Somos-5 sequences}~\cite{gale}. 
(We thank Richard Stanley for telling us about them.)
These are the sequences $a_1,a_2,\dots$ in which any $6$ consecutive terms 
satisfy this relation:
\begin{equation}
\label{eq:somos-5}
a_n a_{n+5} = a_{n+1} a_{n+4} + a_{n+2} a_{n+3} \,. 
\end{equation}
Each term of a Somos-5 sequence is obviously a sub\-tract\-ion-free
rational expression in the first 5 terms $a_1,\dots,a_5\,$. 
More can be shown by extending the arguments in~\cite{gale, malouf}:  
$a_n$ is actually a Laurent polynomial in $a_1,\dots,a_5\,$.
This property is truly remarkable, given the nature of the
recurrence, and the fact that, as $n$ grows, these Laurent polynomials
become huge sums of monomials involving large coefficients;
still, each of these sums cancels out from the
denominator of the recurrence relation
$a_{n+5} = (a_{n+1} a_{n+4} + a_{n+2} a_{n+3})/a_n\,$.
 
We suggest the following analogue of Conjecture~\ref{conj:Laurent}.
  
\begin{conjecture}
\label{conj:Somos-5} 
Every term of a Somos-5 sequence 
is a Laurent polynomial
with nonnegative integer coefficients in the first $5$ terms of the
sequence. 
\end{conjecture}

\section*{Factorization schemes 
}
\label{sec:factorization}

According to Theorem~\ref{th:double diagram TP-criterion}, 
every double wiring diagram 
gives rise to an ``optimal'' total positivity criterion.
We will now show that double wiring diagrams can be used 
to obtain a 
family of bijective parametrizations of the set $G_{>0}$ of all
totally positive matrices; 
this family will include as a special case the parametrization in 
Theorem~\ref{th:imax parametrization}. 

We encode a double wiring diagram by the word of length $n(n-1)$
in the alphabet $\{1, \dots, n-1, \overline{1}, \dots,
\overline{n-1}\}$ obtained by recording the heights of intersections
of pseudolines of like color (traced left to right).
For example, the diagram in Figure~\ref{fig:double-wiring} is encoded
by the word $\overline 2~~ 1~2~~\overline 1~~\overline 2~~1$. 

The words that encode double wiring diagrams have an alternative
description in terms of \emph{reduced expressions} in the symmetric
group~$\mathcal{S}_n\,$. 
Recall that by a famous theorem of E.~H.~Moore \cite{moore},
$\mathcal{S}_n$ is a \emph{Coxeter group} of
type~$A_{n-1}$, i.e., it is generated by the involutions $s_1,\dots,s_{n-1}$
(adjacent transpositions) subject to the relations 
$s_i s_j = s_j s_i$ for $|i-j|\geq 2$, and $s_i s_j s_i = s_j s_i s_i$ for $|i-j| = 1$.
A \emph{reduced word} for a permutation $w\in\mathcal{S}_n$
is a word $\jj=(j_1,\dots,j_l)$ of the shortest possible length~$l=\l(w)$
that satisfies $w=s_{j_1}\cdots s_{j_l}\,$.   
The number $\l(w)$ is called the \emph{length} of~$w$ (it is 
the number of inversions in~$w$). 
The group $\mathcal{S}_n\,$ has a unique element $\wnot$  of maximal length:
the order-reversing permutation of $1, \dots, n$.

It is straightforward to verify that the 
encodings of double wiring diagrams are precisely the shuffles 
of two reduced words for $\wnot\,$, in the barred and unbarred entries
respectively; 
equivalently, these are the reduced words for the
element $(\wnot, \wnot)$ of the Coxeter group~${\mathcal S}_n \times
{\mathcal S}_n\,$. 


\begin{definition}
\label{def:fact-scheme}
{\rm 
A word $\ii$ in the alphabet $\mathcal{A}$  (see (\ref{eq:alphabet})) 
is called a \emph{factorization scheme} if it 
contains each circled entry $\circi$ exactly once,
and the remaining entries encode 
the heights of intersections in a double wiring diagram.

Equivalently, a factorization scheme $\ii$  is a shuffle of two
reduced words for $\wnot$ (one barred and one unbarred) and 
an arbitrary permutation of the entries
$\circled{\,1}\,,\dots,\circled{n}\,$. 
In particular, $\ii$ consists of $n^2$ entries. 
}\end{definition}

To illustrate, the word 
$\ii=\overline 2~~ 1\circled{\,3}\,\,2~~\overline 1\,
\circled{\,1}\,\,\overline 2~~1\circled{\,2}$ 
appearing in Figure~\ref{fig:planar-network} is a factorization
scheme.  

An important example of a factorization scheme 
is the word $\ii_{\max}$ introduced in 
Theorem~\ref{th:imax parametrization}. 
Thus the following result generalizes Theorem~\ref{th:imax parametrization}.

\begin{theorem}
\label{th:TP-bijection}
\cite{FZ}
For an arbitrary factorization scheme
$\ii=(i_1,\dots,i_{n^2})$, 
the product map $x_{\ii}$ given by {\rm (\ref{eq:product-map})}
restricts to a bijection between 
$n^2$-tuples of positive real numbers and 
totally positive $n \times n$ matrices.
\end{theorem} 

\proof
We have already stated that any two double wiring diagrams
are connected by a succession 
of the local ``moves'' shown in Figure~\ref{fig:moves}. 
In the language of factorization schemes, this translates into any two
factorization schemes being connected 
by a sequence  of local transformations of the form
\begin{eqnarray}
\begin{array}{cccl}
\label{eq:iji=jij}
\cdots i\,j\,i \cdots &\leadsto& \cdots j\,i\,j\cdots &,\quad
|i-j|=1\ ,\\[.1in]
\cdots \overline i\,\overline j \,\overline i\cdots &\leadsto&\cdots  
 \overline j\,\overline i\,\overline j\cdots &,\quad |i-j|=1\ ,
\end{array}
\end{eqnarray}
or of the form
\begin{equation}
\label{eq:ab=ba}
\cdots a\,b \cdots \ \leadsto \  \cdots b\,a \cdots \ , 
\end{equation}
where $(a,b)$ is any pair of symbols in~$\mathcal{A}$ different from 
$(i,i\pm 1)$ or $(\overline{i},\overline{i\pm 1})$. 
(This statement is a special case of Tits' theorem~\cite{tits},
for the Coxeter
group~$\mathcal{S}_n\times\mathcal{S}_n\times(\mathcal{S}_2)^n$.) 

In view of Theorem~\ref{th:imax parametrization}, it suffices to show
that if Theorem~\ref{th:TP-bijection} holds for some factorization
scheme $\ii$, then it also holds for any factorization scheme $\ii'$
obtained from $\ii$ by one of the transformations
(\ref{eq:iji=jij})--(\ref{eq:ab=ba}).  
To see this, it is enough to demonstrate that the collections of
parameters $\{t_k\}$ and $\{t'_k\}$ in the equality 
\[
x_{i_1}(t_1)\cdots x_{i_{n^2}}(t_{n^2})  = x_{i'_1}(t'_1)\cdots x_{i'_{n^2}}(t'_{n^2})
\]
are related to each other by (invertible) subtraction-free rational
transformations. 
The latter is a direct consequence of the commutation relations
between elementary  
Jacobi matrices, which can be found in \cite[Section 2.2 and (4.17)]{FZ}.
The most important of these relations are the following.
 
First, for $i=1,\dots,n-1$ and $j = i+1$, we have
\[
x_i (t_1)\,
x_{\circi}(t_2)\,
x_{\circj}(t_3)\,
x_{\overline i} (t_4)
 = 
x_{\overline i}(t'_1) \,
x_{\circi}(t'_2)\,
x_{\circj}(t'_3)  \,
x_i(t'_4) \,, 
\]
where 
\[
t'_1\!=\!\frac{t_3t_4}{T}\,,\ 
t'_2\!=\!{T}\,,\ 
t'_3\!=\!\frac{t_2t_3}{T}\,,\ 
t'_4\!=\!\frac{t_1t_3}{T}\,,\ 
T\!=\!t_2+t_1t_3t_4\,.
\]
The proof of this relation (which is the only nontrivial relation
associated with~(\ref{eq:ab=ba})) 
amounts to verifying that 
\[
\matbr{1}{t_1}{0}{1}\! \matbr{t_2}{0}{0}{t_3}\! \matbr{1}{0}{t_4}{1} 
=
\matbr{1}{0}{t'_1}{1}\!\matbr{t'_2}{0}{0}{t'_3}\!\matbr{1}{t'_4}{0}{1}
. 
\]

Also, for any $i$ and $j$ such that $|i-j|=1$, we have the following
relation associated with (\ref{eq:iji=jij}): 
\begin{equation*}
x_i (t_1) x_j (t_2) x_i (t_3)
= x_j (t'_1) x_i (t'_2) x_j(t'_3)\,, 
\end{equation*}
\begin{equation*}
x_{\overline i} (t_1) x_{\overline j} (t_2) x_{\overline i} (t_3)
= x_{\overline j} (t'_1) x_{\overline i} (t'_2) x_{\overline j}(t'_3)\,, 
\end{equation*}
where
\[
t'_1=\frac{t_2 t_3}{T} \,,\ \ 
t'_2=T \,,\ \ 
t'_3=\frac{t_1 t_2}{T} \,,\ \ 
T=t_1+t_3 \, . 
\]
One sees that in the commutation relations above,
the formulas expressing the $t'_k$ in terms of
the $t_l$ are indeed subtraction-free.
\endproof

Theorem~\ref{th:TP-bijection} suggests an alternative approach to total positivity criteria
via the following \emph{factorization problem}: for a given factorization scheme $\ii$,
find the genericity conditions on a matrix $x$ assuring that $x$ can be factored as
\begin{equation}
\label{eq:x=n2}
x = x_{\ii}(t_1,\dots,t_{n^2})=x_{i_1}(t_1)\cdots x_{i_{n^2}}(t_{n^2}) \,,
\end{equation} 
and compute explicitly the factorization parameters $t_k$ as functions of $x$. 
Then the total positivity of $x$ will be equivalent to the positivity of
all these functions. 
Note that the criterion in Theorem~\ref{th:GP-criterion} was
essentially obtained in this way: 
for the factorization scheme $\ii_{\max}$, the factorization
parameters $t_k$ are Laurent monomials in the initial minors of $x$
(cf.\ Lemma~\ref{lem:initial minors}). 

A complete solution of the factorization problem for an arbitrary factorization scheme
was given in \cite[Theorems~1.9 and 4.9]{FZ}. 
An interesting (and unexpected) feature of this solution is that in
general, the $t_k$ are not Laurent monomials in the minors of~$x$;
the word $\ii_{\max}$ is quite exceptional in this respect. 
It turns out, however, that the $t_k$ are Laurent monomials 
in the minors of another 
matrix $x'$ obtained from $x$ by the following birational transformation:
\begin{equation}
\label{eq:twist}
x'=[x^T\wnot]_+ \wnot (x^T)^{-1} \wnot [\wnot x^T]_-\, .
\end{equation}
Here $x^T$ denotes the transpose of~$x$, and 
$\wnot$ is the permutation matrix with 1's on the antidiagonal; 
finally, $y = [y]_- [y]_0 [y]_+$ denotes the Gaussian (LDU)
decomposition of a square matrix~$y$ provided such a decomposition
exists. 

In the special cases $n=2$ and $n=3$, the transformation $x \mapsto x'$ is
given by
\[
 x'=  
\matbr{x_{11}x_{12}^{-1}x_{21}^{-1}}{x_{21}^{-1}}{x_{12}^{-1}}
{x_{22}\,\det(x)^{-1}}
\]
and
\[
x'=  
\left[ 
{\begin{array}{ccc}
{\displaystyle \frac {{{x}_{1  1}}}{{{x}_{3  1}}\,{{x}_{1  3}}}}
 &  {\displaystyle \frac {\Delta_{12,13}}{{{x}_{3  1}}\,{\Delta_{12,23}}}} & 
{\displaystyle \frac {1}{{{x}_{3  1}}}} \\ [2ex]
{\displaystyle \frac {\Delta_{13,12}}{{\,{{x}_{1  3}}\,\Delta_{23,12}}}} &  
{\displaystyle \frac {
x_{33}\Delta_{12,12}-\det(x)
}{
{\Delta_{23,12}}\,{\Delta_{12,23}}}} & {\displaystyle \frac {{{x}_{3  2}}}{
{\Delta_{23,12}}}} \\ [2ex]
{\displaystyle \frac {1}{{{x}_{1  3}}}} & {\displaystyle \frac {
{{x}_{2  3}}}{{\Delta_{12,23}}}} &   {\displaystyle \frac 
{\Delta_{23,23}}{\det(x)}}
\end{array}}
 \right]  \,. 
\]

The following theorem provides an alternative explanation 
for the family of total positivity criteria in Theorem~\ref{th:double diagram TP-criterion}.

\begin{theorem}
\label{th:factorization-w0}
\cite{FZ}
The right-hand side
of {\rm (\ref{eq:twist})} is well defined for any $x\in G_{>0}\,$;
moreover, the ``twist map'' $x\mapsto x'$ restricts to a bijection 
of $G_{> 0}$ with itself.

Let~$x$ be a totally positive $n\times n$ matrix,
and $\ii$ a factorization scheme. 
Then the parameters $t_1,\dots,t_{n^2}$ appearing in 
{\rm (\ref{eq:x=n2})}
are related by an invertible monomial transformation to the $n^2$
chamber minors (for the double wiring diagram associated with~$\ii$) 
of the twisted matrix~$x'$ given by~{\rm (\ref{eq:twist})}. 
\end{theorem}

In \cite{FZ}, we explicitly describe the monomial transformation in 
Theorem~\ref{th:factorization-w0}, as
well as its inverse, in terms of the combinatorics of the double
wiring diagram.


\section*{Double Bruhat cells}
\label{sec:double Bruhat}

Our presentation in this section will be a bit sketchy;
details can be found in~\cite{FZ}. 

Theorem~\ref{th:TP-bijection} provides a family of bijective (and
biregular) parametrizations 
of the totally positive variety $G_{> 0}$ by $n^2$-tuples of positive
real numbers. 
The totally nonnegative variety $G_{\geq 0}$ is much more complicated 
(note that the map in Theorem~\ref{th:Loewner imax} is surjective but
not injective). 
In this section, we show that $G_{\geq 0}$ splits naturally into 
``simple pieces" corresponding to pairs of permutations
from~$\mathcal{S}_n\,$.


\begin{theorem}
\label{th:uv}
\cite{FZ}
Let $x \in G_{\geq 0}$ be a totally nonnegative matrix.
Suppose that a word $\ii$ in the alphabet $\mathcal{A}$ is such that
$x$ can be factored as $x = x_{\ii} (t_1, \dots, t_m)$ with positive 
$t_1, \dots, t_m$, and $\ii$ has the smallest number of uncircled
entries among all words with this property. 
Then the subword of $\ii$ formed by entries from $\{\overline 1,
\dots, \overline {n-1}\}$ 
(resp.\ from $\{1, \dots, n-1\}$)
is a reduced word for some permutation $u$ (resp.~$v$) 
in~$\mathcal{S}_n$. 
Furthermore, the pair $(u,v)$ is uniquely determined by $x$, i.e.,
does not depend on the choice of~$\ii$.
\end{theorem}

In the situation of Theorem~\ref{th:uv}, we say that $x$ is \emph{of
type}~$(u,v)$. 
Let $G^{u,v}_{> 0} \subset G_{\geq 0}$ denote the subset of all
totally nonnegative matrices of type $(u,v)$; thus $G_{\geq 0}$ is
the disjoint union of these subsets. 

Every subvariety $G^{u,v}_{> 0}$ has a family of parametrizations
similar to those in Theorem~\ref{th:TP-bijection}. 
Generalizing Definition~\ref{def:fact-scheme}, 
let us call a word $\ii$ in the alphabet $\mathcal{A}$ a
\emph{factorization scheme of type}~$(u,v)$ if it contains each
circled entry $\circi$ exactly once, and the barred (resp.\ unbarred)
entries of $\ii$ form a reduced word for~$u$ (resp.~$v$); in
particular, $\ii$ is of length $\l(u)+\l(v)+n$. 

\begin{theorem} 
\label{th:fact-uv}
\cite{FZ}
For an arbitrary factorization scheme $\ii$ of type $(u,v)$, 
the product map $x_{\ii}$ 
restricts to a bijection between 
$(\l(u)+\l(v)+n)$-tuples of positive real numbers and 
totally nonnegative matrices of type $(u,v)$.
\end{theorem}

Comparing Theorems~\ref{th:fact-uv} and \ref{th:TP-bijection},
we see that 
\begin{equation}
\label{eq:TP double cell}
G^{\wnot, \wnot}_{> 0} = G_{> 0} \,,
\end{equation}
i.e., the totally positive matrices are exactly the totally
nonnegative matrices of type~$(\wnot, \wnot)$. 

We now show that the splitting of $G_{\geq 0}$ into the union of
varieties $G^{u,v}_{> 0}$ 
is closely related to the well-known \emph{Bruhat
decompositions} of the general linear group $G = GL_n$.
Let $B$ (resp. $B_-$) denote the subgroup of upper-triangular 
(resp.\ lower-triangular) matrices
in~$G$. 
Recall (see, e.g., \cite[\S4]{alperin-bell})
that each of the double coset spaces $B \backslash G /B$ 
and $B_- \backslash G /B_-$
has cardinality $n!$, and one can choose the permutation matrices $w
\in {\mathcal S}_n$  as their common representatives.
To every two permutations $u$ and $v$ we associate the \emph{double
  Bruhat  cell} 
$G^{u,v} = B u B \cap B_- v B_-$; thus $G$ is the disjoint union of
the double Bruhat cells.

Each set $G^{u,v}$ can be described by equations and
inequalities of the form $\Delta(x)=0$ and/or $\Delta(x)\neq 0$, for
some collection of minors~$\Delta$. 
(See \cite[Proposition~4.1]{FZ} or~\cite{FZcells}.) 
In particular, the open double Bruhat cell $G^{\wnot, \wnot}$ is given by 
non-vanishing of all ``antiprincipal'' minors $\Delta_{[1,i], [n-i+1,n]} (x)$
and $\Delta_{[n-i+1,n], [1,i]} (x)$ for $i = 1, \dots, n\!-\!1$. 

\begin{theorem}
\label{th:TNN-bijection}
\cite{FZ}
A totally nonnegative matrix is of type $(u,v)$ 
if and only if it belongs to the double Bruhat cell~$G^{u,v}$.  
\end{theorem}


In view of (\ref{eq:TP double cell}), Theorem~\ref{th:TNN-bijection} 
provides the following 
simple test for total positivity of a totally
nonnegative matrix.

\begin{corollary}
\label{cor:TP-w0}
\cite{gasca-pena}
A totally nonnegative matrix~$x$ is totally positive 
if and only if $\Delta_{[1,i], [n-i+1,n]} (x) \neq 0$
and $\Delta_{[n-i+1,n], [1,i]} (x)\neq 0$ for $i = 1, \dots, n$.
\end{corollary}

The results 
obtained above for $G^{\wnot, \wnot}_{> 0} = G_{> 0}$ 
(as well as their proofs) extend 
to the variety $G^{u,v}_{> 0}$ for
an arbitrary pair of permutations $u,v\in\mathcal{S}_n\,$.
In particular, the factorization schemes for $(u,v)$ (or rather 
their uncircled parts) can be visualised by  
\emph{double wiring diagrams of type} $(u,v)$ in the same way as before, 
except now any two pseudolines intersect \emph{at most} once, and 
the lines are permuted ``according to $u$ and~$v$.'' 
Every such diagram has $\l(u)+\l(v)+n$ chamber minors, and their positivity 
provides a criterion for a matrix $x \in G^{u,v}$ to belong
to~$G^{u,v}_{> 0}\,$. 
The factorization problem and its solution provided by
Theorem~\ref{th:factorization-w0} 
extend to any double Bruhat cell, with an appropriate modification of the 
twist map $x \mapsto x'$.
The details can be found in~\cite{FZ}.

If the double Bruhat cell containing a matrix $x \in G$ is not
specified,  
then testing $x$ for total nonnegativity becomes a much harder problem;
in fact,  
every known criterion involves exponentially many (in~$n$) minors. 
(See \cite{CLR} for related complexity results.) 
The following corollary of a result by Cryer~\cite{cryer76} was given
by Gasca and Pe\~na \cite{gasca-pena-siam}. 


\begin{theorem}
\label{th:GP-TNN criterion}
An invertible square matrix is totally nonnegative if and only if all
its minors occupying 
several initial rows or 
several initial columns are nonnegative,
and all its leading principal minors are positive. 
\end{theorem}

This criterion involves $2^{n+1}-n-2$ minors, which is roughly the
square root of the total number of minors. 
We do not know whether this criterion is optimal.  

\section*{Oscillatory matrices}

We conclude the paper by discussing the intermediate class of
oscillatory matrices 
that was introduced and intensively studied by
Gantmacher and Krein~\cite{GK35, GK}. 
A matrix is \emph{oscillatory} 
if it is totally nonnegative while some power of it is
totally positive; thus the set of oscillatory matrices 
contains $G_{> 0}$ and is contained in~$G_{\geq 0}\,$. 
The following theorem provides several equivalent characterizations
of oscillatory matrices; the equivalence of (a)-(c) was proved in 
\cite{GK}, while the rest of the conditions were given in \cite{FZosc}. 

\begin{theorem} 
\label{th:GK classical}
{\rm \cite{GK, FZosc}}
For an invertible totally nonnegative $n\times n$ matrix~$x$,
the following are equivalent: 
\begin{itemize}
\item[(a)] 
$x$ is oscillatory; 
\item[(b)]
$x_{i,i+1} > 0$ and $x_{i+1,i} > 0$ for $i = 1, \dots, n\!-\!1$; 
\item[(c)]
$x^{n-1}$ is totally positive;
\item[(d)] 
$x$ is \emph{not} block-triangular (cf.\ Figure~\ref{fig:block-triangular}); 
\item[(e)]
$x$ can be factored as
$x=x_\ii(t_1,\dots,t_l)$, for positive $t_1,\dots,t_l$ 
and a word $\ii$ that contains every symbol of the form~$i$ or
$\bar i$ at least once;
\item[(f)]
$x$ lies in a double Bruhat cell $\,G^{u,v}$, where both $u$ and $v$
do not fix any set $\{1,\dots,i\}$, for $i=1,\dots,n\!-\!1$. 
\end{itemize} 
\end{theorem}


\begin{figure}[ht]
\setlength{\unitlength}{1pt} 
\begin{center}
$
\left[
\begin{array}{ccccc}
* & * & 0 & 0 & 0  \\
* & * & 0 & 0 & 0  \\
* & * & * & * & *  \\
* & * & * & * & *  \\
* & * & * & * & *  \\
\end{array}
\right]
\quad
\left[
\begin{array}{cccccc}
* & * & * & * & *  \\
* & * & * & * & *  \\
0 & 0 & * & * & *  \\
0 & 0 & * & * & *  \\
0 & 0 & * & * & *  \\
\end{array}
\right]
$
\end{center}
\vspace{-0.1in}
\caption{
Block-triangular matrices
}
\label{fig:block-triangular}
\end{figure}

\vspace{-.1in}

\proof
Obviously, 
${\rm (c)}\Longrightarrow{\rm (a)}\Longrightarrow{\rm (d)}$. 
Let us prove the equivalence of (b), (d), and~(e). 
By Theorem~\ref{th:Loewner classical},
$x$ can be represented as the weight
matrix of some planar network $\Gamma(\ii)$ with positive edge
weights. Then (b) means that sink $i+1$ (resp.~$i$) can be reached from
source~$i$ (resp.\ $i+1$), for all~$i$; 
(d) means that for any $i$,
at least one sink $j>i$ is reachable from a source $h \leq i$, 
and at least one sink $h \leq i$ is reachable from a source~$j>i$;
and (e) means that $\Gamma(\ii)$ contains positively- and
negatively-sloped edges connecting any two consecutive levels $i$
and~$i+1$.
These three statements are easily seen to be  equivalent. 

By Theorem~\ref{th:TNN-bijection}, ${\rm (e)}\Longleftrightarrow{\rm (f)}$. 
It remains to show that ${\rm (e)}\Longrightarrow{\rm (c)}$. 
In view of Theorem~\ref{th:fact-uv} and (\ref{eq:TP double cell}),
this can be restated as follows: given any permutation $\jj$ of the
entries $1,\dots,n\!-\!1$, 
prove that the concatenation $\jj^{n-1}$ of $n\!-\!1$ copies
of~$\jj$ contains a reduced word for~$\wnot\,$.
Let $\jj'$ denote the subsequence of $\jj^{n-1}$ constructed
as follows. 
First, $\jj'$ contains all $n\!-\!1$ entries of $\jj^{n-1}$ 
which are equal to~$n\!-\!1$.  
Second, $\jj'$ contains the $n\!-\!2$ entries equal to $n\!-\!2$
which interlace the $n\!-\!1$ entries chosen at the previous step. 
We then include $n\!-\!3$ interlacing entries equal to $n\!-\!3$, etc. 
The resulting word $\jj'$ of length $\binom{n}{2}$ 
will be a reduced word for~$\wnot\,$,
for it will be equivalent, under the transformations
(\ref{eq:ab=ba}), to the lexicographically maximal reduced word
$\jj_{\max}=(n\!-\!1,n\!-\!2,n\!-\!1, n\!-\!3,n\!-\!2,n\!-\!1, \dots)$. 
\endproof



\section*{Acknowledgments}

We thank Sara Billey for suggesting a number of editorial
improvements. 
This work was supported in part 
by NSF grants \#DMS-9625511 and \#DMS-9700927.

\end{document}